\documentstyle[amscd,amssymb,verbatim,diagrams,12pt]{amsart}
\newcommand{\CYBE}{\operatorname{CYBE}}
\newcommand{\reg}{{\operatorname{reg}}}
\newcommand{\res}{{\operatorname{res}}}
\newcommand{\const}{{\operatorname{const}}}
\newcommand{\fg}{{\mathfrak g}}
\newcommand{\fn}{{\mathfrak n}}
\newcommand{\fb}{{\mathfrak b}}
\newcommand{\fh}{{\mathfrak h}}

\newcommand{\Perf}{\operatorname{Perf}}
\newcommand{\ev}{{\operatorname{ev}}}

\newcommand{\bA}{{\bf A}}
\newcommand{\bB}{{\bf B}}

\newcommand{\bM}{{\bf M}}
\newcommand{\bN}{{\bf N}}
\newcommand{\bP}{{\bf P}}

\newcommand{\bH}{{\bf H}}

\newcommand{\ZZ}{{\cal Z}}

\newcommand{\ssl}{\mathfrak{sl}}

\renewcommand{\mod}{\operatorname{mod}}

\newcommand{\und}{\underline}

\newcommand{\OO}{{\cal O}}
\newcommand{\RR}{{\cal R}}

\newcommand{\II}{{\cal I}}

\newcommand{\BB}{{\cal B}}

\newcommand{\tors}{\operatorname{tors}}

\newcommand{\G}{{\Bbb G}}

\newcommand{\mg}{{\frak m}}
\newcommand{\hra}{\hookrightarrow}
\newcommand{\lan}{\langle}
\newcommand{\ran}{\rangle}

\newcommand{\CC}{{\cal C}}

\newcommand{\tr}{\operatorname{tr}}

\newcommand{\Mat}{\operatorname{Mat}}

\newcommand{\Proj}{\operatorname{Proj}}
\renewcommand{\P}{{\Bbb P}}

\newcommand{\si}{\sigma}

\newcommand{\de}{\delta}
\newcommand{\eps}{\epsilon}

\numberwithin{equation}{subsection}

\newtheorem{thm}{Theorem}[subsection]
\newtheorem{prop}[thm]{Proposition}
\newtheorem{lem}[thm]{Lemma}
\newtheorem{cor}[thm]{Corollary}
{  \theoremstyle{definition}
\newtheorem{defi}[thm]{Definition}
\newtheorem{ex}[thm]{Example}

\newtheorem{rem}[thm]{Remark}

}

\newcommand{\Pf}{\noindent {\it Proof}}
\newcommand{\id}{\operatorname{id}}

\renewcommand{\AA}{{\cal A}}

\newcommand{\FF}{{\cal F}}

\newcommand{\VV}{{\cal V}}

\newcommand{\LL}{{\cal L}}

\newcommand{\Om}{\Omega}

\newcommand{\Hom}{\operatorname{Hom}}

\newcommand{\Ext}{\operatorname{Ext}}
\newcommand{\End}{\operatorname{End}}
\newcommand{\Res}{\operatorname{Res}}

\newcommand{\Aut}{\operatorname{Aut}}

\renewcommand{\a}{\alpha}
\renewcommand{\b}{\beta}
\newcommand{\om}{\omega}
\newcommand{\De}{\Delta}
\newcommand{\la}{\lambda}
\newcommand{\th}{\theta}
\newcommand{\C}{{\Bbb C}}

\newcommand{\Z}{{\Bbb Z}}

\newcommand{\Ga}{\Gamma}
\newcommand{\e}{\operatorname{e}}
\newcommand{\wt}{\widetilde}
\newcommand{\ot}{\otimes}

\newcommand{\sub}{\subset}
\newcommand{\ed}{\qed\vspace{3mm}}

\title{Geometrization of trigonometric solutions of the associative and classical Yang-Baxter equations}
\author{Alexander Polishchuk}
\address{
    Department of Mathematics, 
    University of Oregon, 
    Eugene, OR 97403, USA; HSE University, Russian Federation; and Korea Institute for 
    Advanced Study 
  }
  \email{apolish@@uoregon.edu}
\thanks{Supported in part by the NSF grant and by the Russian Academic Excellence Project `5-100' within the framework of the HSE University Basic Research Program}

\begin{document}

\begin{abstract}
We describe a geometric construction of all nondegenerate trigonometric solutions of the associative and classical Yang-Baxter equations.
In the associative case the solutions come from symmetric spherical orders over the irreducible nodal curve of arithmetic genus $1$, while
in the Lie case they come from spherical sheaves of Lie algebras over the same curve.
\end{abstract}

\maketitle

\section*{Introduction}

Recall that the {\it classical Yang-Baxter equation (CYBE)} for a Lie algebra $\fg$ is the equation
\begin{equation}\label{CYBE-intro-eq}
[r^{12}(v),r^{13}(v+v')]+[r^{12}(v),r^{23}(v')]+[r^{13}(v+v'),r^{23}(v')]=0
\end{equation}
on a meromorphic function $r(v)$ in a neighborhood of zero, taking values in $\fg\ot \fg$,
where $r^{12}=r\ot 1\in U(\fg)^{\ot 3}$, etc.
This is a well studied equation related to the theory of classical integrable systems and of quantum groups
(see e.g., \cite{Drinfeld}, \cite{ES}). It is usually coupled with the unitarity condition
$$r^{21}(-v)=-r(v).$$
Belavin and Drinfeld \cite{BD} showed that in the case when $\fg$ is simple, 
all nondegenerate solutions of the CYBE are either elliptic or
trigonometric or rational and classified elliptic and trigonometric solutions.

The {\it associative Yang-Baxter equation (AYBE)} for an associative algebra $A$ is the equation 
\begin{equation}\label{AYBE-eq}
r^{12}(-u',v)r^{13}(u+u',v+v')-
r^{23}(u+u',v')r^{12}(u,v)+
r^{13}(u,v+v')r^{23}(u',v')=0,
\end{equation}
where $r(u,v)$ is a meromorphic function of two complex variables in a neighborhood of $(0,0)$, taking values in $A\ot A$,
where $r^{12}=r\ot 1\in A^{\ot 3}$, etc. It is usually coupled with the skew-symmetry condition
$$r^{21}(-u,-v)=-r(u,v).$$
In the above form the AYBE was introduced in \cite{pol02}; the constant version was introduced
in \cite{Aguiar}. In \cite{P-nodal-massey} we proved an analog of Belavin-Drinfeld classifications for nondegenerate
skew-symmetric solutions of the AYBE for the matrix algebra $\Mat_n(\C)$ in terms of some combinatorial data,
called {\it associative Belavin-Drinfeld data (BD data)}
(for the definition of the nondegeneracy condition, which is stronger than the one
for the CYBE, see \cite[Def.\ 1.4.3]{LP-YB}).

We are interested in geometric constructions of solutions of \eqref{CYBE-intro-eq} and \eqref{AYBE-eq}.
Our starting point is the construction (going back to  \cite{pol02}) of solutions of the AYBE for $\Mat_n(\C)$
(resp., of the CYBE for the Lie algebra $\ssl_n$) coming from a pair of families of $1$-spherical objects. 
It was shown in \cite{pol02} that all nondegenerate elliptic solutions arise in this way from the families $(\VV_s)$ and $(\OO_x)$
on an elliptic curve, where $(\VV_s)$ are stable bundles of given rank and degree.
The natural problem is to construct geometrically all trigonometric solutions.

In \cite{P-nodal-massey} we showed that some of trigonometric solutions of the AYBE are realized geometrically
using families of $1$-spherical objects $(\VV_s)$, $(\OO_x)$ on some nodal Calabi-Yau curves, where $\VV_s$ are simple vector bundles.
In \cite{LP-YB} we realized all nondegenerate trigonometric solutions of the AYBE using objects in the Fukaya categories
of square-tiled surfaces.
 
In this work we present an algebro-geometric realization of all nondegenerate trigonometric solutions of both the AYBE and the CYBE,
using appropriate sheaves over the irreducible nodal curve of arithmetic genus $1$.

For the AYBE we use the framework of (symmetric) spherical orders over a projective integral curves developed in \cite{P-ainf-orders}. Essentially, such an order is a coherent sheaf of $\OO$-algebras $\AA$ over a curve $C$, with $h^0(C,\AA)=1$, generically isomorphic to a matrix algebra and equipped with a (symmetric) perfect pairing
$$\AA\ot \AA\to \om_C,$$
where $\om_C$ is the dualizing sheaf (see Sec.\ \ref{sph-orders-sec} for details).
In the category of $\AA$-modules we have pairs of $1$-spherical objects $(\AA\ot\LL, V\ot \OO_x)$, where $\LL$ is a line bundle over $C$, $x$ is a smooth point of $C$,
$V$ a vector space such that $\AA|_x\simeq \End(V)$. Therefore, one gets the corresponding solution of the AYBE.
Our first main result (see Theorem \ref{trig-AYBE-thm}) shows that the trigonometric solution of the AYBE associated with associative BD data 
comes in this way from a naturally constructed symmetric spherical order over the projective nodal curve of arithmetic genus $1$.

For the CYBE there is a well known framework of acyclic sheaves of Lie algebras discovered by Cherednik \cite{Cherednik} and
developed by Burban-Galinat \cite{BG}. The fact that all trigonometric solutions arise in this way is mentioned in \cite{Drinfeld},
however, it seems that aside from some examples, no general construction of the corresponding sheaves of Lie algebras over $C$ existed before
this work. By analogy with the associative operad case, we introduce the notion of a {\it (symmetric) spherical sheaf of Lie algebras} as a coherent sheaf of Lie algebras $\LL$ over $C$,
with $H^*(C,\LL)=0$, equipped with a (symmetric) perfect pairing
$$\LL\ot \LL\to \om_C.$$
As in the constructions of \cite{Cherednik} and \cite{BG}, near every smooth point of $C$, we get the corresponding Manin triple, or equivalently, a classical $r$-matrix
(see Sec.\ \ref{spherical-Lie-sec}).
Our second main result is the construction of a symmetric spherical sheaf of Lie algebras giving rise to a given 
nondegenerate trigonometric solution of the CYBE for a simple Lie algebra
$\fg$ (see Theorem \ref{main-Lie-thm}). Note that it complements nicely the result of Burban-Galinat \cite{BG} that all nondegerate {\it rational} solutions of the CYBE come from
appropriate sheaves of Lie algebras over the cuspidal cubic (the geometric construction of elliptic solutions goes back to \cite{Cherednik}). Finally, we should mention that our construction of a spherical sheaf of Lie algebras from a solution of the CYBE
is vaguely reminiscent of the construction of algebro-geometric spectral data from algebras of commuting differential operators, dating back to Mumford's paper \cite{M}.

The paper consists of two parts: in Section 1 we consider the associative structures (spherical orders and the AYBE), while Section 2 is devoted to Lie structures 
(spherical sheaves of Lie algebras and the CYBE). 
After recalling in Sec.\ \ref{1-sph-solutions-sec} and \ref{sph-orders-sec} how solutions of the AYBE appear from spherical orders, in Sec.\ \ref{sph-isotropic-sec}
we give a general construction of spherical orders on the irreducible nodal curve $C$ of arithmetic genus 1 from certain maximal isotropic
subalgebras in $\Mat_n(k)\times \Mat_n(k)$. Then in Sec. \ref{isotr-parab-sec} we give a construction of such isotropic subalgebras starting from a pair of
parabolic subalgebras in $\Mat_n(k)$ equipped with the isomorphism between their semisimple quotients.
Then, after reminding in Sec.\ \ref{BD-data-sec} the definition of associative BD data and the classification of trigonometric solutions of the AYBE,
in Sec.\ \ref{BD-sph-order-sec} we construct a spherical order on $C$ starting from BD data. In Sec.\ \ref{AYBE-comp-sec} we compute
the corresponding solution of the AYBE (see Theorem \ref{trig-AYBE-thm}).
In Sec.\ \ref{more-sph-sec} we prove some structure results on symmetric spherical orders. In particular, we prove that every such order on the irreducible
nodal curve of arithmetic genus $1$ comes from some maximal isotropic subalgebra of $\Mat_n(k)\times \Mat_n(k)$ as in Sec.\ \ref{sph-isotropic-sec}.

In Sec.\ \ref{formal-CYBE-sec} we review some standard results on the CYBE and (infinite-dimensional) Manin triples.
In Sec.\ \ref{spherical-Lie-sec} we introduce the notion of a symmetric spherical sheaf of Lie algebra and discuss its relation to Manin triples.
After reminding the Belavin-Drinfeld's classification of trigonometric solutions of the CYBE in Sec.\ \ref{BD-class-sec},
in Sec.\ \ref{sph-Lie-nodal-sec} we show that every such solution comes from a spherical sheaf of Lie algebras on the irreducible nodal curve of
arithmetic genus $1$ (see Theorem \ref{main-Lie-thm}).

\medskip

\noindent
{\it Conventions}. Our ground field $k$ is assumed to be algebraically closed of characteristic zero. In the parts dealing with classification
of trigonometric solutions we assume $k=\C$.

\medskip

\noindent
{\it Acknowledgments}. I am grateful to Riley Casper for a useful discussion of orders on nodal curves and to the anonymous referee
for many useful comments.
Part of this work was done during a visit to the Hebrew University of Jerusalem, which I would like to thank for hospitality.


\section{Spherical orders and the AYBE}

\subsection{Solutions of the AYBE associated with two families of $1$-spherical objects}\label{1-sph-solutions-sec}

Let us recall the general construction of solutions of the AYBE (\cite{pol02}, \cite{LP-YB}) from $1$-Calabi-Yau $A_\infty$-categories.
Assume we are given such a minimal $A_\infty$-category $\CC$ and two sets of isomorphism classes
of objects in $\CC$, $\mathcal{X}$ and $\mathcal{Y}$,
such that for every pair of distinct objects $x_1,x_2 \in \mathcal{X}$ (resp. $y_1,y_2 \in \mathcal{Y}$), $\Hom^*(x_1,x_2)=0$ (resp. $\Hom^*(y_1,y_2)=0$). We also assume 
that $\Hom^{\neq 0}(x,y)=0$ (and so $\Hom^{\neq 1}(y,x)=0$) for $x\in \mathcal{X}$ and $y\in \mathcal{Y}$. Then 
dualizing the triple product 
$$m_3:\Hom^0(x_2,y_2) \ot\Hom^1(y_1,x_2)\ot\Hom^0(x_1,y_1)\to\Hom^0(x_1,y_2),$$
where $x_1,x_2\in \mathcal{X}$, $y_1,y_2\in \mathcal{Y}$, using the Calabi-Yau pairing, we get a tensor 
$$r^{x_1,x_2}_{y_1,y_2}:\Hom^0(x_2,y_2)\ot \Hom^0(x_1,y_1)\to\Hom^0(x_1,y_2)\ot \Hom^0(x_2,y_1).$$
Using the $A_\infty$-relations and the cyclic symmetry, one can check that this tensor satisfies
the {\it general AYBE},
\begin{equation}\label{Rid0}
(r^{x_1x_2}_{y_1y_3})^{13}(r^{x_2x_3}_{y_2y_3})^{12}+
(r^{x_3x_2}_{y_1y_2})^{23}(r^{x_1x_3}_{y_1y_3})^{13}-
(r^{x_1x_3}_{y_2y_3})^{12}(r^{x_1x_2}_{y_1y_2})^{23}=0,
\end{equation}
viewed as an equation between maps
\[ \Hom^0(x_3,y_3) \ot \Hom^0(x_2,y_2) \ot \Hom^0(x_1,y_1) \to  \Hom^0(x_2,y_3) \ot \Hom^0(x_1,y_2) \ot \Hom^0(x_3,y_1), \]
and the following {\it skew-symmetry condition}:
\begin{equation}\label{skew-sym-eq}
(r^{x_1x_2}_{y_1y_2})^{21}=-r^{x_2x_1}_{y_2y_1}.
\end{equation}

Choosing identifications of all spaces $\Hom^*(x_i,y_j)$ with the same vector space $V$ allows us to view $r^{x_1,x_2}_{y_1,y_2}$
as an element of $\End(V)\ot \End(V)$. If in addition, the parameter spaces ${\mathcal X}$ and ${\mathcal Y}$ are subsets of $\C$ and
the $r$ depends only on the differences $x_1-x_2$, $y_1-y_2$, we get a solution of the AYBE \eqref{AYBE-eq}.
More precisely, to get \eqref{AYBE-eq}, one has to swap the second and the third factors in the triple tensor product and use
the skew-symmetry condition (see \cite[Introduction]{LP-YB}).

In examples leading to nondegenerate solutions the objects in $\mathcal{X}$ and $\mathcal{Y}$ are in addition {\it $1$-spherical} in the sense of
\cite{ST} (see \cite[Sec.\ 1.5]{pol02} for the relevance of this assumption).

\subsection{Solutions of the AYBE associated with spherical orders}\label{sph-orders-sec}

Let $C$ be an integral projective curve over $k$. We denote by $\om_C$ the dualizing sheaf on $C$.
By an {\it order} over $C$ we mean a torsion free
coherent sheaf of $\OO_C$-algebras $\AA$ such that $\AA_{\eta}$ is a matrix algebra over the field of
rational functions on $C$ (where $\eta$ is the generic point of $C$).
Let us recall the following definitions from \cite{P-ainf-orders} (see \cite[Def.\ 0.1.3]{P-ainf-orders}
and \cite[Prop.\ 3.2.2]{P-ainf-orders}).

\begin{defi}\label{sph-order-def}
An order $\AA$ over $C$, such that $H^0(C,\OO_C)=k$, is called 
{\it spherical} if $\AA$ is a $1$-spherical object in the perfect derived category of left $\AA$-modules, $\Perf(\AA)$, 
or equivalently, for a nonzero morphism (unique up to rescaling) $\tau:\AA\to \om_C$ of $\OO_C$-modules,
the induced morphism
$$\nu:\AA\to\und{\Hom}(\AA,\om_C):y\mapsto (x\mapsto \tau(xy))$$
is an isomorphism (equivalently, one can consider the morphism $y\mapsto (x\mapsto \tau(xy))$).
A spherical order is called {\it symmetric} if in addition one has $\tau(xy)=\tau(yx)$.
\end{defi}

It is proved in \cite[Prop.\ 3.2.2]{P-ainf-orders} that an order $\AA$ is spherical if and only if $\AA^{op}$ is spherical.
Note that by Serre duality, for a spherical order $\AA$ one has 
$$h^1(\AA)=h^1(\und{\Hom}(\AA,\om_C))=h^0(\AA),$$
so $\chi(\AA)=0$.

\begin{lem}\label{trace-order-lem} 
Let $\AA$ be a symmetric spherical order over $C$ and let $p\in C$ be a smooth point such that
$\AA|_p\simeq \End(V)$, where $V$ is a vector space. Then $\tau|_p:\AA|_p\to \om_C|_p$ can be identified
with a nonzero multiple of the trace map, $\tr:\End(V)\to k$, i.e., there exists a nonzero element $\a_p\in \om_C|_p$
such that $\tau|_p(?)=\tr(?)\cdot \a_p$.
\end{lem}

\Pf . This follows from \cite[Prop.\ 3.2.2(ii)]{P-ainf-orders}.
\ed

Suppose we are 
given an order $\AA$ over an integral curve $C$
and a line bundle $M$ on $C$ such that $H^*(C,\AA\ot M)=0$. Then for any smooth points $x\neq x'$ in $C$ we 
define a map
\begin{equation}\label{pre-r-matrix-map}
\rho_\AA(x,x';M): \AA|_x\ot M|_x\ot T_xC\to \AA|_{x'}\ot M|_{x'}
\end{equation}
from the commutative diagram
\begin{diagram} 
\AA|_x\ot M|_x\ot T_xC&\lTo{\res_x}& H^0(\AA\ot M(x))\\
&\rdTo{\rho_\AA(x,x';M)}&\dTo_{\ev_{x'}}\\
&&\AA|_{x'}\ot M|_{x'}
\end{diagram}
where the map $\res_x$ is given by the restriction to $x$:
$$H^0(\AA\ot M(x))\to \AA|_x\ot M(x)|_x\simeq \AA|_x\ot M|_x\ot T_xC,$$
where we use the canonical isomorphism $\OO_C(x)|_x\simeq T_xC$ dual to the residue map.
Note that the map $\res_x$ an isomorphism due to the assumption $H^*(C,\AA\ot M)$=0.
We can view $\rho_\AA(x,x';M)$ as an element of 
$$(\AA^\vee\ot\om_C)|_x\ot \AA|_{x'}\ot M|_{x'}\ot M^{-1}|_x\simeq \AA|_{x'}\ot\und{\Hom}(\AA,\om_C)|_x\ot M|_{x'}\ot M^{-1}|_x.$$
In particular, applying $\nu^{-1}|_x$ we get an element
$$(\id\ot \nu^{-1}|_x)\rho_\AA(x,x';M)\in \AA|_{x'}\ot \AA|_x \ot M|_{x'}\ot M^{-1}|_x.$$

Now assume that $\AA$ is a symmetric spherical order, and we are given
a smooth open subset $U\sub C$, together with a trivialization
\begin{equation}\label{triv-A-U-eq}
\AA|_U\simeq\End(V)\ot \OO_U,
\end{equation}
and in addition let $\LL$ be a family of line bundles of degree $0$ on $C$, with the base $S$,
such that for generic members $L=\LL_s$, $L'=\LL_{s'}$ of this family one has 
$h^0(\AA\ot L'\ot L^{-1})=0$.
Note that for any point $x\in U$ we can view $V\ot \OO_x$ as a right module over $\AA$,
using the standard left $\End(V)$-module structure on $V$ and the trivialization \eqref{triv-A-U-eq}.
Also, by Lemma \ref{trace-order-lem}, there exists a uniquely defined nowhere vanishing $1$-form $\a\in \om_C(U)$,
such that under the trivialization \eqref{triv-A-U-eq}, the morphism $\tau|_U$ gets identified with
$$\End(V)\ot \OO_U\to \om_U: A\mapsto \tr(A)\cdot \a.$$ 
Hence, the isomorphism \eqref{triv-A-U-eq} is compatible with duality in the form of the commutative diagram
\begin{equation}\label{nu-a-compatibility}
\begin{diagram}
\AA|_x&\rTo{\nu|_x}&(\AA^\vee\ot \om_C)|_x\\
\dTo{\sim}&&\dTo{\sim}\\
\End(V)&\rTo{\a_x}&\End(V)^*\ot \om_C|_x
\end{diagram}
\end{equation}
where in the bottom horizontal arrow we use the standard selfduality of $\End(V)$.

The Serre duality on the perfect derived category of $\AA$-modules is discussed in \cite[Sec.\ 3.2]{P-ainf-orders}.
In particular, in the above situation the $\AA$-modules $\AA\ot \LL_s$ and $V\ot\OO_x$ are $1$-spherical.
For example, for a perfect $\AA$-module $P$, the composition and $\tau:\AA\to \om_C$ induce a canonical pairing
\begin{equation}\label{order-Serre-pairing}
\Ext^1_{\AA}(P,\AA\ot \LL_s)\ot \Hom_{\AA}(\AA\ot\LL_s,P)\to\Ext^1_{\AA}(\AA\ot \LL_s,\AA\ot\LL_s)\simeq 
H^1(\AA)\rTo{\tau} H^1(\om_C)\simeq k.
\end{equation}

Thus, we can apply the construction of solutions of the AYBE from Sec.\ \ref{1-sph-solutions-sec} to the families of objects 
\begin{equation}\label{two-families-eq}
(\AA\ot \LL_s)_{s\in S}, \ \ (V\ot \OO_x)_{x\in U}.
\end{equation}
Note that the minimal $A_\infty$-structure on the corresponding subcategory $\Perf(\AA)$ obtained by
homological perturbation can be chosen to be cyclic with respect to the above canonical pairings induced by the Serre duality. 
Indeed, this is proved in the same way as 
in the proof of \cite[Cor.\ C]{P-ainf-orders}.

For a line bundle $L=\LL_s$ and a point $x\in U$, we have an identification
\begin{equation}\label{Hom-order-V-eq}
\Hom_{\AA}(\AA\ot L, V\ot \OO_x)\simeq L^{-1}|_x\ot V.
\end{equation}
On the other hand, using Serre duality \eqref{order-Serre-pairing},
we get an identification
\begin{equation}\label{Ext-order-V-eq}
\Ext^1_{\AA}(V\ot \OO_x,\AA\ot L)\simeq \Hom_{\AA}(\AA\ot L,V\ot \OO_{x})^*\simeq L|_x\ot V^*.
\end{equation}

\begin{lem}\label{serre-dual-lem}
For $x\in U$, let us consider the composition
$$\Hom_{\AA}(V\ot \OO_x,\AA\ot L(x)|_x\ot \OO_x)\rTo{\de} \Ext^1_{\AA}(V\ot \OO_x,\AA\ot L)\rTo{\sim} L|_x\ot V^*,$$
where the second arrow is the isomorphism \eqref{Ext-order-V-eq} and $\de$ is the connecting homomorphism associated
with the exact sequence
\begin{equation}\label{AL(x)-ex-seq}
0\to \AA\ot L\to \AA\ot L(x)\to \AA\ot L(x)|_x\ot \OO_x\to 0.
\end{equation}
Then under the isomorphism \eqref{triv-A-U-eq}
this composition gets identified with the map
$$\Hom_{\AA}(V\ot \OO_x,\AA\ot L(x)|_x\ot \OO_x)\simeq\Hom_{\End(V)}(V,\End(V))\ot L|_x\ot T_C|_x\rTo{\id\ot \a_x} V^*\ot L|_x.$$
\end{lem}

\Pf . To simplify the notation, we assume that $L$ is trivial; the proof in the general case is essentially the same.
The Serre duality isomorphism is given by the composition
$$\Ext^1_{\AA}(V\ot\OO_x,\AA)\to\Hom_{\AA}(\AA,V\ot\OO_x)^*\ot\Ext^1_{\AA}(\AA,\AA)\to \Hom_{\AA}(\AA,V\ot\OO_x)^*,$$
where the first arrow is the dualization of the composition map, while the second arrow is induced by the
isomorphism
$$\Ext^1_{\AA}(\AA,\AA)\simeq H^1(\AA)\rTo{\tau} H^1(\om_C)\simeq k.$$
Now we use the commutative diagram
\begin{equation}\label{serre-conn-hom-comp-diagram}
\begin{diagram}
\Hom_{\AA}(V\ot \OO_x,\AA(x)|_x\ot\OO_x)&\rTo{}&\Hom_{\AA}(\AA,V\ot\OO_x)^*\ot\Hom_{\AA}(\AA,\AA(x)|_x\ot\OO_x)\\
\dTo{\de}&&\dTo{}\\
\Ext^1_{\AA}(V\ot\OO_x,\AA)&\rTo{}&\Hom_{\AA}(\AA,V\ot\OO_x)^*\ot\Ext^1_{\AA}(\AA,\AA)
\end{diagram}
\end{equation}
where the vertical arrows are given by the connecting homomorphisms associated with the exact sequence
\eqref{AL(x)-ex-seq}.
Note that the top horizontal arrow can be identified with the natural map
$$\Hom_{\End(V)}(V,\End(V))\ot T_C|_x\hra V^*\ot \End(V)\ot T_C|_x\rTo{\sim} V^*\ot \AA(x)|_x.$$
Since the composition
$$V^*\simeq \Hom_{\End(V)}(V,\End(V))\hra V^*\ot \End(V)\rTo{\id\ot\tr} V^*$$
is equal to the identity, by the commutativity of diagram \eqref{serre-conn-hom-comp-diagram}, it remains to check that the composed map
$$\AA(x)|_x\to H^1(\AA)\rTo{\tau} H^1(\om_C)\simeq k$$
is given by $\AA(x)|_x\simeq \End(V)\to T_C|_x\rTo{\tr\ot\a_x} k$.
But this follows immediately from the definition of $\a$ and from the commutativity of the diagram
\begin{diagram}
\AA(x)|_x&\rTo{}& H^1(\AA)\\
\dTo{\tau}&&\dTo{\tau}\\
\om_C(x)&\rTo{}& H^1(\om_C)
\end{diagram}
\ed

\begin{prop}\label{AYBE-order-sol-prop} 
Let $\AA$ be a symmetric spherical order over an integral curve $C$, 
The solution of the AYBE associated with the families \eqref{two-families-eq} and the trivialization
\eqref{triv-A-U-eq},
is given by
$$r^{L,L'}_{x,x'}=(\id\ot\nu^{-1}|_x)\rho_\AA(x,x';L'\ot L^{-1})\in \End(V)\ot\End(V)\ot L'|_{x'}\ot L^{-1}|_{x'}\ot (L')^{-1}|_x\ot L|_x,$$
where $x\neq x'$ and $(L,L')$ are such that $H^*(\AA\ot L'\ot L^{-1})=0$.
\end{prop}

\Pf . For brevity let us write $\AA\ot L=\AA L$, $V\ot \OO_x=V\OO_x$, etc.
Recall that $r^{L,L'}_{x,x'}$ is obtained by dualizing the triple product
$$m_3:\Hom_{\AA}(\AA L',V\OO_{x'})\ot \Ext^1_{\AA}(V\OO_x,\AA L')\ot \Hom_{\AA}(\AA L,V\OO_x)\to\Hom_{\AA}(\AA L,V\OO_{x'})$$
and using  identifications \eqref{Hom-order-V-eq}, \eqref{Ext-order-V-eq}.

\medskip

\noindent
{\bf Step 1. Interpretation as a Massey product}.
The above triple product is given by the well-defined and univalued triple Massey product $MP$,
so using the standard recipe for its calculation we have to include the canonical morphism
$$\Ext^1_{\AA}(V\OO_x,\AA L')\ot V\OO_x\to \AA L'[1]$$
into an exact triangle 
\begin{equation}\label{Massey-pr-ex-tr}
\AA L'\to C \to \Ext^1\ot V\OO_x \to \AA L'[1],
\end{equation}
where we abbreviate $\Ext^1_{\AA}(V\OO_x,\AA L')$ as $\Ext^1$. Then we use the fact that the maps
$$\Hom_{\AA}(\AA L,C)\rTo{r_x} \Hom_{\AA}(\AA L,\Ext^1\ot V\OO_x)\simeq \Ext^1\ot \Hom_{\AA}(\AA L,V\OO_x),$$
$$\Hom_{\AA}(C,V\OO_{x'})\rTo{i_{x'}} \Hom_{\AA}(\AA L',V\OO_{x'})$$
are isomorphisms and compute $m_3=MP$ as
\begin{align*}
& \Hom_{\AA}(\AA L',V\OO_{x'})\ot \Ext^1\ot \Hom_{\AA}(\AA L,V\OO_x)\rTo{i_{x'}^{-1}\ot r_x^{-1}}
 \Hom_{\AA}(C,V\OO_{x'})\ot \Hom_{\AA}(\AA L,C) \\ 
& \rTo{\kappa} \Hom_{\AA}(\AA L,V\OO_{x'}),
\end{align*}
where $\kappa$ is given by the composition.
The map $\kappa$ is compatible with restrictions to $x'$ and the composition map
$$\kappa_{x'}: \Hom_{\AA|_{x'}}(C|_{x'},V)\ot \Hom_{\AA|_{x'}}(\AA|_{x'} L|_{x'},C|_{x'})\to \Hom_{\AA|_{x'}}(\AA|_{x'}L|_{x'},V).$$
On the other hand, the isomorphism $i_{x'}$ is compatible with the isomorphism $\AA|_{x'} L'|_{x'}\rTo{\sim} C|_{x'}$ induced
by the triangle \eqref{Massey-pr-ex-tr}, which we still denote as $i_{x'}$.
It follows that $m_3=MP$ can be rewritten as the following composition
\begin{align*} 
&((L')^{-1}_{x'}\ot V)\ot \Ext^1\ot \Hom_{\AA}(\AA L,V\OO_x)\rTo{\id \ot r_x^{-1}}
((L')^{-1}_{x'}\ot V)\ot \Hom_{\AA}(\AA L,C)\rTo{\id \ot i_{x'}^{-1}}\\
& ((L')^{-1}_{x'}\ot V)\ot \Hom_{\AA|_{x'}}(\AA|_{x'} L|_{x'},\AA|_{x'}L'|_{x'})\simeq ((L')^{-1}_{x'}\ot V)\ot (L^{-1}|_{x'}L'|_{x'}\ot \AA|_{x'}) \\
&\to L^{-1}|_{x'}\ot V,
\end{align*}
where the last arrow is induced by the isomorphism $\AA|_{x'}\simeq \End(V)$ coming from \eqref{triv-A-U-eq}.
In other words, the dualization of $m_3$ corresponds to the composition
\begin{align*}
&\Ext^1\ot \Hom_{\AA}(\AA L,V\OO_x))\rTo{r_x^{-1}}\Hom_{\AA}(\AA L,C)\to \Hom_{\AA|_{x'}}(\AA|_{x'} L|_{x'},C|_{x'})
\rTo{i_{x'}^{-1}}\\ 
&L^{-1}|_{x'}\ot L'|_{x'}\ot \AA|_{x'}\simeq L^{-1}|_{x'}\ot L'|_{x'}\ot \End(V).
\end{align*}

\medskip

\noindent
{\bf Step 2. Identification of the exact triangle}.
We claim that using the identification \eqref{Ext-order-V-eq}, the exact triangle \eqref{Massey-pr-ex-tr}
can be identified with an exact triangle of the form
\begin{equation}\label{AL'x-triangle}
\AA L'\rTo{\iota} \AA L'(x)\rTo{r} L'|_x\ot V^*\ot V\OO_x\to \AA L'[1],
\end{equation}
where $\iota$ is the natural inclusion, and the map $r$ comes from the isomorphism 
\begin{equation}\label{A(x)-EndV-identification}
\AA (x)|_x\simeq \End(V)\ot T_C|_x\rTo{\a_x} V^*\ot V.
\end{equation}
Indeed, we define \eqref{AL'x-triangle} as the triangle isomorphic to the standard exact triangle
$$\AA L'\rTo{\iota} \AA L'(x) \to \AA L'(x)|_x\ot \OO_x \to \AA L'[1]$$
using the isomorphism \eqref{A(x)-EndV-identification}.

Let us set for brevity $A:=V\OO_x$, $B:=\AA L'$, $W:=L'|_x\ot V^*$, $C'=\AA L'(x)$.
We have constructed an exact triangle of the form
$$B\to C'\to W\ot A\to B[1]$$
The map $W\ot A\to B[1]$ corresponds to some map $\b:W\to \Ext^1(A,B)$, and we have a morphism of exact triangles
\begin{diagram}
B&\rTo{}& C'&\rTo{}& W\ot A&\rTo{}& B[1]\\
\dTo{\id}&&\dTo{}&&\dTo{\b\ot\id}&&\dTo{\id}\\
B&\rTo{}& C&\rTo{}& \Ext^1(A,B)\ot A&\rTo{}& B[1]
\end{diagram}
Note that we check that $\b$ is an isomorphism then it would follow that this is an isomorphism of exact triangles.
Thus, it remains to compute $\b$. More precisely, we need to show that 
$\b$ coincides with the isomorphism
$$L'|_x\ot V^*\rTo{\sim} \Ext^1_{\AA}(V \OO_x,\AA L')$$
given by \eqref{serre-dual-lem}. Indeed, $\b$ is obtained by applying $\Hom_{\AA}(V\OO_x,?)$
to a morphism of degree $1$ in the exact triangle \eqref{AL'x-triangle}. Note that this triangle
comes from the exact sequence isomorphic to \eqref{AL(x)-ex-seq}, with $L'$ instead of $L$,
via the isomorphism $\AA L'(x)|_x\simeq L'|_x\ot V^*\ot V$ coming from 
\eqref{A(x)-EndV-identification}. Hence, $\b$ is equal to the composition 
\begin{align*}
&L'|_x\ot V^*\rTo{\sim}\Hom_{\AA}(V\OO_x,L'|_x\ot V^*\ot V\OO_x)\rTo{\sim}
\Hom_{\AA}(V\OO_x,\AA L'(x)|_x \OO_x)\rTo{\de}\\
&\Ext^1_{\AA}(V\OO_x,\AA L'), 
\end{align*}
where the second arrow is induced by \eqref{A(x)-EndV-identification} and the last arrow
is the connecting homomorphism for the exact sequence. Now the claim follows from 
Lemma \ref{serre-dual-lem}.

\medskip

\noindent
{\bf Step 3. Conclusion of the computation}.
Using the isomorphism of exact triangle from Step 2, we get the following commutative diagram
\begin{diagram}
\Ext^1\ot\Hom_{\AA}(\AA L,V\OO_x)&\lTo{r_x}&\Hom_{\AA}(\AA L,C)&\rTo{i_{x'}^{-1}}& \AA|_{x'}\ot (L^{-1}L')|_{x'}\\
\uTo{\sim}&&\uTo{\sim}&&\uTo{=}\\
\AA|_x\ot (L^{-1}L')|_x&\lTo{\a_x\cdot \res_x}&\Hom_{\AA}(\AA L,\AA L'(x))&\rTo{\ev_{x'}}&\AA|_{x'}\ot (L^{-1}L')|_{x'}
\end{diagram}
where the left vertical arrow comes from the identifications $\AA|_x\simeq \End(V)$, 
$$\Ext^1\ot\Hom_{\AA}(\AA L,V\OO_x)\simeq L'|_x\ot V^*\ot L^{-1}|_x\ot V.$$
It remains to observe that after inverting $r_x^{-1}$, the composition of obtained arrows in the top row gives
the dualization of $m_3$ (by Step 1). On the other hand, by commutativity of the diagram
\eqref{nu-a-compatibility}, the similar composition in the bottom row corresponds to $(\id\ot\nu^{-1}|_x)\rho_\AA(x,x';L'\ot L^{-1})$,
as asserted.
\ed



\subsection{Spherical orders over the irreducible nodal curve of arithmetic genus $1$ associated with maximal isotropic subalgebras}\label{sph-isotropic-sec}

Let $C$ be the irreducible nodal curve of arithmetic genus $1$, and let $\pi:\P^1\to C$ denote the normalization map
such that $q=\pi(0)=\pi(\infty)$ is the node on $C$. We denote by $(x_0:x_1)$ the homogeneous coordinates on $\P^1$,
where $x_1=0$ at the point $0\in \P^1$ and $x_0=\infty$ at the point $\infty\in\P^1$.

Let us equip  $\Mat_n(k)\times\Mat_n(k)$ with the symmetric pairing given by 
\begin{equation}\label{tr-diff-pairing}
\lan (a,b), (a',b')\ran=\tr(aa')-\tr(bb').
\end{equation}

\begin{defi} Let $\VV$ be a vector bundle of rank $n$ over $\P^1$, and
let $I\sub \End(\VV|_0)\oplus \End(\VV|_\infty)$ be a maximal isotropic subalgebra with respect to the pairing 
\eqref{tr-diff-pairing}.
We denote by 
$$\AA(\VV,I)\sub \pi_*\und{\End}(\VV)$$ 
the subsheaf consisting of sections $a$ such that 
$(a(0),a(\infty))\in I$.
\end{defi}

\begin{lem}\label{sym-sph-order-isotr-lem} Let $\AA=\AA(\VV,I)$. Assume that $h^0(\AA)=1$. Then $\AA$ is a symmetric
spherical order over $C$. Furthermore, the center of $\AA$ is $\OO_C$.
\end{lem}

\Pf . Since $(1,1)\in I$, for any $(a,b)\in I$ we have 
$$0=\lan (a,b),(1,1)\ran=\tr(a)-\tr(b).$$
It follows that the natural map $\tr:\und{\End}(\VV)\to \pi_*\OO_{\P^1}$ restricts to a map
$$\tr:\AA\to \OO_C\simeq \om_C.$$
It satisfes $\tr(xy)=\tr(yx)$, so we just have to check that the induced morphism
$\AA\to\und{\Hom}(\AA,\OO_C)$ is an isomorphism. This is clear away from the node, so
it is enough to study the situation in the formal neighborhood of the node.
Let 
$$R=\hat{\OO}_{C,q}\sub k[[x]]\oplus k[[y]]$$ denote the completed local ring of $C$ at the node, and let 
$A=\hat{\AA}_q$ be the completion of $\AA$ at the node.
By definition, $A$ is a subalgebra in $\Mat_n(k)[[x]]\oplus\Mat_n(k)[[y]]$ consisting of $(a(x),b(y))$ such that
$(a(0),b(0))\in I$. We just have to check that given $(a(x),b(y))\in \Mat_n(k)((x))\oplus\Mat_n(k)((y))$,
such that for any $(a'(x),b'(y))\in A$ one has $(\tr(aa'),\tr(bb'))\in R$, we necessarily have $(a,b)\in A$.
Indeed, first, using the inclusion
$$x\Mat_n(k)[[x]]\oplus y\Mat_n(k)[[y]]\sub A,$$ we deduce that
$a(x)\in \Mat_n(k)[[x]]$ and $b\in \Mat_n(k)[[y]]$.
Next, for any $(a',b')\in I\sub A$, the condition $(\tr(a(x)a'),\tr(b(y)b'))\in R$ is equivalent to
$$\tr(a(0)a')=\tr(b(0)b')),$$
i.e., that $(a(0),b(0))$ is orthogonal to $(a',b')$. Since $I$ is maximal isotropic, this implies that $(a(0),b(0))\in I$.

For the last assertion, we first observe that there is an inclusion $\OO_C\sub \ZZ$, where $\ZZ$ denotes the center of $\AA$.
Since over $U=C\setminus q$, the sheaf $\AA$ coincides with 
$\pi_*\und{\End}(\VV)$, we have $\ZZ|_U=\OO_U$. Hence, the quotient $\ZZ/\OO_C$ is a torsion sheaf.
On the other hand, we claim that $\AA\cap \pi_*\OO_{\P^1}=\OO_C$.
Indeed, this amounts to checking that $(\la,\mu)\in I$, with $\la,\mu\in k$, only if $\la=\mu$. But this immediately follows
from the condition 
$$\tr(\la)-\tr(\mu)=n(\la-\mu)=0.$$

Thus, we have an embedding 
$$\AA/\OO_{C}\sub \pi_*\und{\End}(\VV)/\pi_*\OO_{\P^1}.$$
Since the latter sheaf is torsion free, it follows that $\AA/\OO_C$ is torsion free.
But $\ZZ/\OO_C$ is a subsheaf of $\AA/\OO_C$, hence, $\ZZ/\OO_C=0$.
\ed

\begin{rem}\label{1-0-rem} 
We needed the fact that $(1,0)\not\in I$ in Lemma \ref{sym-sph-order-isotr-lem} to prove that $\AA=\OO_C$.
In positive characteristic dividing $n$ this is not necessarily true.
Indeed, there are examples of maximal isotropic subalgebras
in $\Mat_n(k)$ and we can just take $I$ to be a direct sum of two such subalgebras. 
\end{rem}

\begin{rem} In Sec.\ \ref{more-sph-sec} we will show that conversely, every symmetric spherical order on $C$
is of the form $\AA(\VV,I)$, with $\VV$ of specific form.
\end{rem}

\subsection{Maximal isotropic subalgebras associated with parabolic subalgebras}\label{isotr-parab-sec}

To construct isotropic subalgebras $I$ we will use parabolic subalgebras in $\Mat_n(k)$, i.e.,
subalgebras conjugate to a standard parabolic subalgebra of the form
$$\bP_+(\Pi)=\bB_+ +\sum_{\a\in \Pi} \lan e_{-\a}\ran,$$
where $\bB_+$ is the subalgebra of upper-triangular matrices and
$\Pi$ is a subset of positive weights, closed under sums and such that if $\a+\b\in\Pi$ with $\a,\b>0$
then $\a,\b\in\Pi$. Here we identify nonzero weights with pairs $\a=(i,j)$, $i\neq j$, so that $\a>0$ (resp., $\a<0$) if and only if
$i<j$ (resp., $i>j$), and $-\a=(j,i)$.
We will also consider parabolic subalgebras of the form
$$\bP_-(\Pi)=\bB_- +\sum_{\a\in\Pi} \lan e_{\a}\ran,$$
where $\bB_-$ is the subalgebra of lower-triangular matrices.

For a parabolic subalgebra $\bP$, let $J(\bP)$ be its Jacobson radical, 
$L(\bP):=\bP/J(\bP)$, the semisimple quotient of $\bP$, which is a product of matrix algebras, and we
denote by $\pi_L:\bP\to L(\bP)$ the natural projection.

Note that for a standard parabolic subalgebra we have
$$L(\bP_+(\Pi))=\bH+\sum_{\a\in\Pi}(\lan e_\a\ran+\lan e_{-\a}\ran),$$
where $\bH$ is the subalgebra of diagonal matrices. In particular,
$$\dim \bP_+(\Pi)=\frac{n(n+1)}{2}+|\Pi|, \ \ \dim L(\bP_+(\Pi))=n+2|\Pi|.$$

\begin{defi}
Let $(\bP_1,\bP_2)$ be a pair of parabolic subalgebras, such that we have an algebra isomorphism
$$\phi:L(\bP_1)\rTo{\sim} L(\bP_2),$$
Then we define a subalgebra $I(\bP_1,\bP_2,\phi)\sub\Mat_n(k)\oplus\Mat_n(k)$ as follows:
$$I(\bP_1,\bP_2,\phi):=\{(a,b)\in \bP_1\oplus \bP_2 \ |\ \pi_L(b)=\phi(\pi_L(a))\}.$$
\end{defi}

\begin{lem}\label{parabolic-max-isotr-lem}
$I(\bP_1,\bP_2,\phi)$ is a maximal isotropic subalgebra with respect to the pairing \eqref{tr-diff-pairing}.
\end{lem}

\Pf . Indeed, it is straightforward to see that $I(\bP_1,\bP_2,\phi)$ is a subalgebra.
Thus, to check that it is isotropic, it is enough to check that $\tr(a)=\tr(b)$ for any $(a,b)\in I(\bP_1,\bP_2,\phi)$. 
Since $\tr(a)=\tr(\pi_L(a))$, this follows from the fact that $\phi$ is 
compatible with traces, i.e.,
$$\tr(\phi(x))=\tr(x).$$
But both $L(\bP_1)$ and $L(\bP_2)$ can be identified with the same standard Levi subalgebra $L\sub \Mat_n(k)$
conisting of block diagonal matrices. It remains to note that every automorphism of $L$ is inner, so it preserves the trace.

To calculate the dimension of $I(\bP_1,\bP_2,\phi)$, note that for any parabolic subalgebra $\bP$, we have
$$\dim \bP=\frac{n^2+\dim L(\bP)}{2}$$
(this can be checked for standard parabolic subalgebras).
Hence, $\dim \bP_1=\dim \bP_2$, and 
$$\dim I(\bP_1,\bP_2,\phi)=\dim \bP_1+\dim \bP_2-\dim L(\bP_1)=2\dim \bP_1-\dim L(\bP_1)=n^2.$$
\ed

\begin{prop}\label{parabolic-constr-prop} 
For a vector bundle $\VV$ of rank $n$ over $\P^1$, a pair of parabolic subalgebras
$\bP_1\sub\End(\VV|_0)$, $\bP_2\sub\End(\VV|_\infty)$, and 
an algebra isomorphism $\phi:L(\bP_1)\rTo{\sim} L(\bP_2)$, let $\AA=\AA(\VV,I(\bP_1,\bP_2,\phi))$ be the corresponding
order over $C$. If $h^0(\AA)=1$ then $\AA$ is a symmetric spherical order over $C$ with the center $\OO_C$.
\end{prop}

\Pf . By Lemma \ref{parabolic-max-isotr-lem}, $I=I(\bP_1,\bP_2,\phi)$ is a maximal isotropic subalgebra.
Now the result follows from 
Lemma \ref{sym-sph-order-isotr-lem}.
\ed

\begin{rem}
The result of Proposition \ref{parabolic-constr-prop} still holds in positive characteristic.
Indeed, we have $(1,0)\not\in I$, since $\pi_L(1)=1\neq 0$ (see Remark \ref{1-0-rem}). 
\end{rem}

\subsection{Associative BD data and trigonometric solutions of the AYBE}\label{BD-data-sec}

Recall that in \cite{P-nodal-massey} we classified trigonometric nondegenerate solutions of the AYBE
in terms of {\it associative BD data} $(\si_0,\si,\Ga_1,\Ga_2)$, where $\si_0,\si\in S_n$ is a pair of cyclic permutations
and
$$\Ga_1,\Ga_2\sub \Ga_{\si_0}=\{(i,\si_0(i))\in [1,n]^2 \ | \ i\in [1,n]\}$$
are proper subsets such that $\si(\Ga_1)=\Ga_2$, where we let permutations act on $[1,n]^2$ diagonally.
These data are considered up to a permutation, so we will always assume that $\si_0$ is standard: $\si_0(i)=i+1$, so
$$\Ga_{\si_0}=\{(1,2),(2,3)\ldots,(n,1)\}.$$
We will usually suppress $\si_0$ from the notation and just write the associative BD data as $(\si,\Ga_1,\Ga_2)$.

To describe the formula for the corresponding associative $r$-matrix,
let us extend the subsets $\Ga_1,\Ga_2$ to the subsets 
$\Pi_1,\Pi_2\sub [1,n]^2$, where for $a=1,2$,
$$\Pi_a=\{(i,\si_0^r(i)) \ |\ (i,\si_0(i))\in \Ga_a, (\si_0(i),\si_0^2(i))\in\Ga_a,\ldots,(\si_0^{r-1}(i),\si_0^r(i))\in\Ga_a\}.$$

The main property of $\Pi_1$ and $\Pi_2$ is that whenever $i,j,k$ are in cyclic order (with respect to $\si_0$) then
$$(i,k)\in \Pi_a \iff (i,j)\in \Pi_a \land (j,k)\in \Pi_a,$$
for  $a=1,2$.
It follows that we still have
\begin{equation}\label{permutation-Pi-12-eq}
\si(\Pi_1)=\Pi_2.
\end{equation}

Recall that solutions of \eqref{AYBE-eq} are obtained from solutions $r^{x_1,x_2}_{y_1,y_2}$ of \eqref{Rid0}, \eqref{skew-sym-eq},
such that $r$ depends only on the differences $x_1-x_2,y_1-y_2$.
For the description of trigonometric solutions, it is convenient to consider the intermediate form of dependence, where 
$r$ depends on $x_1-x_2,y_1,y_2$. Actually, we will use multiplicative variables $x_i$, and will present a formula
for $r(\la;y_1,y_2)$, with $\la\in k^*$, such that $r(x_1/x_2;y_1,y_2)$ satisfies \eqref{Rid0}, and such that
in the case $k=\C$, the solution $r(e^u;e^{v_1},e^{v_2})$ is equivalent to a trigonometric solution that depends only on $u$ and $v_1-v_2$ and appears in the classification theorem \cite[Thm.\ 0.1]{P-nodal-massey}
(this equivalence is proved in \cite[Lem.\ 6.1]{P-nodal-massey}).

We use the standard basis $e_{ij}$ of the matrix algebra $\Mat_n(k)$. We also write $h_i=e_{ii}$.
In the formula below we denote pairs 
$(i,j)$ with $i\neq j$ by the variable $\a$.
We write $(i,j)>0$ (resp., $(i,j)<0$) if $i<j$ (resp., $i>j$).
Also, for $\a=(i,j)$ we set $-\a=(j,i)$.
We set
\begin{equation}\label{r-const-AYBE-sol}
\begin{array}{l}
r_{\si,\const}(\la,z)=\frac{1}{1-z} \sum_{\a>0} e_{-\a}\otimes e_{\a}+
\frac{z}{1-z} \sum_{\a>0} e_{\a}\otimes e_{-\a}+\\
\frac{1}{1-z}\sum_{i=1}^n h_i\otimes h_i-
(1-\la^n)^{-1}\sum_{i=1}^n\sum_{k=0}^{n-1} \la^k h_i\otimes h_{\si^k(i)}
\end{array}
\end{equation}
Now the solution of the AYBE associated with the associative BD data $(\si,\Ga_1,\Ga_2)$ is given by
\footnote{Our formula differs from the \cite[Eq. (36)]{P-nodal-massey} by swapping $x$ and $y$ (and by a sign, which is not
essential). This is due to different conventions
in forming the $r$-matrix from a Massey product. Without this swap the $r$-matrix satisfies a slightly different equation from 
\eqref{Rid0}, see the discussion around Eq. (0-7) in \cite[Introduction]{LP-YB}.}
\begin{equation}\label{r-gen-AYBE-sol}
\begin{array}{l}
r_{\si,\Ga_1,\Ga_2}(\la;x,y)=r_{\si,\const}(\la,x/y)+\\
\sum_{\a>0,k\ge 1}[\la^{-k}e_{-\tau^k(\a)}\otimes e_{\a}-\la^k e_{\a}\otimes e_{-\tau^k(\a)}]+\\
\sum_{\a<0,k\ge 1}[x\la^{-k} e_{-\tau^k(\a)}\otimes e_{\a}-y\la^k e_{\a}\otimes e_{-\tau^k(\a)}].
\end{array}
\end{equation}
Here we use the operation $\tau$ defined only on $\Pi_1\sub [1,n]\times [1,n]$ and given by $\tau(\a)=\si(\a)$;
the summation is extended only over those $(k,\a)$ for which $\tau^k(\a)$ is defined. 

\begin{ex} Consider the case $n=6$, $\si=(136245)$, $\Ga_1=\{(6,1),(1,2)\}$, $\Ga_2=\{(2,3),(3,4)\}$.
Then we have one more element in $\Pi_1$, namely $(6,2)$. The map $\tau$ is given by
$$\tau(6,1)=(2,3), \ \tau(1,2)=(3,4), \ \tau(6,2)=(2,4).$$
There are no $\a$ for which $\tau^2(\a)$ would be defined.
The only positive element of $\Pi_1$ is $(1,2)$, whereas $(6,1)$ and $(6,2)$ are negative.
Hence, in this case we have
\begin{align*}
&r_{\si,\Ga_1,\Ga_2}(\la;x,y)=r_{\si,\const}(\la,x/y)+\la^{-1} e_{43}\ot e_{12}-\la e_{12}\ot e_{43} \\
&+x\la^{-1}(e_{32}\ot e_{61}+e_{42}\ot e_{62})-y\la (e_{61}\ot e_{32}+e_{62}\ot e_{42}).
\end{align*}
\end{ex}

\subsection{From associative BD data to spherical orders over $C$}
\label{BD-sph-order-sec}

\begin{ex}
In the case when $\Ga_i=\emptyset$, for every cyclic permutation $\si\in S_n$, 
we define the subalgebra $I_\si\sub\Mat_n(k)\oplus \Mat_n(k)$ by 
$$I_\si=I(\bB_+,\bB_-,\si)=
\{(a_{ij}),(b_{ij}) \ |\ a_{ij}=0 \text{ for } i>j, \ b_{ij}=0 \text{ for } i<j, \ a_{ii}=b_{\si(i)\si(i)} \text{ for all } i\},$$
and we set 
$$\AA_\si=\AA(\OO_{\P^1}^{\oplus n},I_\si).$$
Note that global sections of $\AA_\si$ can be viewed as matrices $a=(a_{ij})\in\Mat_n(k)$ such that $(a,a)\in I_\si$.
This easily implies that $a$ can only be a scalar matrix.
Hence, by Proposition \ref{parabolic-constr-prop},
$\AA_\si$ is a symmetric spherical order with the center $\OO_C$.
\end{ex}

Now let us consider the general case, so we start with associative BD data $(\si,\Ga_1,\Ga_2)$.
Without loss of generality we can assume that $(n,1)\not\in \Ga_1$ (since $\Ga_1$ is a proper subset of $\Ga_{\si_0}$).
It follows that for every $(i,j)\in \Pi_1$ one has $i<j$.
To construct our order we need one additional choice: we choose $m$ such that $(\si_0^{-1}(m),m)\not\in \Ga_2$.
Note that such $m$ exists since $\Ga_2$ is a proper subset of $\Ga_{\si_0}$.

Now we define the sequence of $0$'s and $1$'s numbered by $i\in [1,n]$:
$$\eps_i:=\begin{cases} 0, & i<m, \\ 1, i\ge m.\end{cases}$$
We are going to construct our spherical order in the form
$\AA=\AA(\VV,I(\bP_1,\bP_2,\si))$ with 
$$\VV=\bigoplus_{i=1}^n \OO(\eps_i),$$
for certain parabolic subalgebras $\bP_1,\bP_2$ in $\Mat_n(k)$ and an isomorphism between their Levi components
induced by the permutation $\si$.
We will use the standard trivializations
of the fibers of $\OO(1)$ at $0$ and $\infty$ (given by the homogeneous coordinates $x_0$ and $x_1$ on $\P^1$, respectively)
to get identifications of $\VV|_0$ and $\VV|_\infty$ with $\Mat_n(k)$.

We set 
$$\bP_1:=\bP_+(\Pi_1)=\{(a_{ij}) \ |\ a_{ij}=0 \text{ for } i>j, (j,i)\not\in\Pi_1\}.$$

Let $S_m\sub [1,n]^2$ denote the set of pairs $(i,j)$
such that either $i<j<m$ or $m\le i<j$ or $i\ge m>j$. In other words, $(i,j)\in S_m$ if and only if
$\si_0^{-m+1}i<\si_0^{-m+1}j$. Note that $S_m$ has the following properties:
\begin{itemize}
\item 
if $(i,j)\in S_m$ and $(j,k)\in S_m$ then $(i,k)\in S_m$;
\item
if $(i,j)\in S_m$ then for any $k$, either $(i,k)\in S_m$  or $(k,j)\in S_m$;
\item
if $(i,j)\not\in S_m$ and $i\neq j$ then $(j,i)\in S_m$.
\end{itemize}

Note that our assumption on $m$ is equivalent to $(n,1)\not\in \si_0^{-m+1}\Pi_2$.
This implies that we have an inclusion
$$\Pi_2\sub S_m.$$ 

Now we set
$$\bP_2=\si_0^{m-1}\bP_-(\si_0^{-m+1}\Pi_2)=
\{(a_{ij}) \ |\ a_{ij}=0 \text{ for } (i,j)\in S_m, (i,j)\not\in\Pi_2\}.$$
Here we view any permutation $\si'$ as an automorphism of $\Mat_n(k)$ such that 
$\si'(e_{ij})=e_{\si'(i),\si'(j)}$.

Finally, we observe that the Levi components of our subalgebras are
$$L(\bP_r)={\mathrm{span}}((h_i \ |\ i=1,\ldots,n), (e_{ij},e_{ji} \ |\ (i,j)\in\Pi_r)), \ \text{ for } r=1,2.$$ 
Hence, due to the property \eqref{permutation-Pi-12-eq}, we have 
$$\si L(\bP_1)=L(\bP_2).$$
Thus, we can form the maximal isotropic subalgebra $I=I(\bP_1,\bP_2,\si)$, and consider the order
$\AA=\AA(\VV,I(\bP_1,\bP_2,\si))$ on $C$.


Explicitly, the subalgebra $I=I(\bP_1,\bP_2,\si)$ consists of pairs of matrices $(a,b)$ satisfying the following conditions
\begin{enumerate}
\item $j<i,(j,i)\not\in \Pi_1\implies a_{ij}=0$, 
\item $(i,j)\in \Pi_1 \text{ or } i=j \implies a_{ij}=b_{\si(i)\si(j)}, a_{ji}=b_{\si(j)\si(i)}$,
\item $(i,j)\in S_m, (i,j)\not\in \Pi_2 \implies b_{ij}=0$.
\end{enumerate}

It will be convenient to use the following description of $\bP_i$.
Let $\bB_1^+=\bB^+$ (resp., $\bB_2^-=\bB^-$) be the subalgebras of upper-triangular (resp., lower-triangular)
matrices, and set 
$$\bB_2^\pm=\si_0^{m-1}\bB^\pm,$$ 
so that $a\in \bB_2^+$ (resp., $a\in \bB_2^-$) if and only if $a_{ij}=0$ for $(j,i)\in S_m$ (resp., for $(i,j)\in S_m$).
Let us also set for $i=1,2$,
$$\bA_i^+=\sum_{\a\in\Pi_i}\lan e_\a\ran, \ \ \bA_i^-=\sum_{\a\in\Pi_i}\lan e_{-\a}\ran.$$
Then we have 
$$\bP_1=\bB_1^+ + \bA_1^-, \ \ \bP_2=\bB_2^-+ \bA_2^+.$$
Note the Levi components are given by
\begin{equation}\label{Levi-P-eq}
L(\bP_r)=\bH + \bA_r^+ + \bA_r^-, \ \text{ for } r=1,2,
\end{equation}
where $\bH\sub \Mat_n(k)$ is the subspace of diagonal matrices.

Also, let $\bN_i^\pm$ denote the unipotent radicals in $\bB_i^\pm$.
We denote by $\pi_{\bB_i^\pm}$, $\pi_{\bN_i^\pm}$, $\pi_{\bA_i^\pm}$ and $\pi_\bH$ the coordinate projectors
to the corresponding subspaces.
Let us also introduce the following subspaces in $\Mat_n(k)$:
\begin{align*}
&\bM(\OO(-1)):=\lan e_{ij} \ |\ i<m,j\ge m\ran, \ \
\bM(\OO):=\lan e_{ij} \ |\ i<m,j<m \text{ or } i\ge m,j\ge m\ran, \\ 
&\bM(\OO(1)):=\lan e_{ij} \ |\ i\ge m>j\ran,
\end{align*}
and let $\pi_{\OO(-1)}$, $\pi_{\OO}$ and $\pi_{\OO(1)}$ be the corresponding 
coordinate projectors on $\Mat_n(k)$. 
Note that all the above projectors commute, and we have
$$\pi_{\bB_1^+}\pi_{\OO(1)}=\pi_{\bB_2^-}\pi_{\OO(1)}=0.$$
Also, we have
$$\pi_{\bB_1^\pm}\pi_{\OO}=\pi_{\bB_2^\pm}\pi_{\OO}.$$

\begin{lem}\label{nilpotent-operator-lem} 
The operators 
$$\th^+:=\si\pi_{\bA_1^+}\pi_\OO, \ \ \th^-:=\si^{-1}\pi_{\bA_2^-}\pi_\OO$$ 
are nilpotent.
\end{lem}

\Pf .  Let us check this for $\th^+$ (for $\th^-$ the argument is the same).
We have $\th^+(e_{ij})=0$ unless $i<j<m$ or $m\le i<j$, and $(i,j)\in\Pi_1$. If these conditions are satisfied then
$$\th^+(e_{ij})=\pi_\OO(e_{\si(i)\si(j)}),$$ 
which is either $0$ or $e_{\si(i)\si(j)}$. 
Since there exists a positive power $k$ such that $(\si^k(i),\si^k(j))\not\in\Pi_1$, the assertion follows.
\ed

Everywhere in the calculations below we will view sections of $\OO(1)$ on $\P^1$ using the trivialization
over $\P^1\setminus\{\infty\}$ given by $x_0$. For example, to a global section $c_0x_0+c_1x_1$ we associate
the linear function $c_0+c_1t$, where $t=x_1/x_0$ is the coordinate on $\P^1\setminus\{\infty\}$.
In terms of this identification the evaluation at $0$ sends $c_0+c_1t$ to $c_0$, while the evaluation at $\infty$
sends $c_0+c_1t$ to $c_1$ (since we use $x_1$ to trivialize $\OO(1)$ at $\infty$).
Later we will similarly consider rational sections of $\OO(1)$, which in terms of our trivialization correspond to
rational functions of degree $\le 1$ in $t$.

\begin{lem}\label{trivial-h0-lem} 
One has $h^0(\AA)=1$. Hence, $\AA$ is a symmetric spherical order over $C$ with the center $\OO_C$.
\end{lem}

\Pf . A global section of $\AA$ is an $n\times n$ matrix $a=(a_{ij})$, with $a_{ij}\in H^0(\OO(\eps_i-\eps_j))$,
such that $(a(0),a(\infty))$ belongs to $I(\bP_1,\bP_2,\si)$.
Thus, we can write $a=a_0+a_1t$, where
$$\pi_{\OO(-1)}a_0=0, \ \ a_1\in\bM(\OO(1)).$$
Note that
$$a(0)=a_0\in \bP_1, \ \ a(\infty)=\pi_{\OO}a_0+a_1\in \bP_2.$$
Furthermore, due to the description \eqref{Levi-P-eq} of the Levi components of $\bP_1$ and $\bP_2$, the condition 
$(a(0),a(\infty))\in I(\bP_1,\bP_2,\si)$ implies that
$$\pi_\bH a(0)=\si^{-1}\pi_\bH a(\infty), \ \ \pi_{\bA_1^\pm}a(0)=\si^{-1}\pi_{\bA_2^\pm}a(\infty).$$

Since $\pi_{\bB_1^+}\pi_{\OO(1)}=0$, $a_0=\pi_{\OO}a_0+\pi_{\OO(1)}a_0$ and
$a_1\in\bM(\OO(1))$, we have
\begin{equation}\label{B1+a0-relation}
\pi_{\bB_1^+}a_0=\pi_{\bB_1^+}\pi_{\OO}a_0=\pi_{\bB_1^+}\pi_{\OO}a(\infty)=\pi_{\bB_2^+}\pi_{\OO}a(\infty)=
\pi_\OO\pi_{\bA_2^+}a(\infty)+\pi_\bH a(\infty).
\end{equation}

Next, let us analyze the condition
\begin{equation}\label{pi-A1-A2-si-eq}
\pi_{\bA_1^-}a_0=\si^{-1}\pi_{\bA_2^-}a(\infty).
\end{equation}
We have $\pi_{\bA_2^-}a_1=0$ since $a_1\in\bM(\OO(1))$, hence, 
$$\pi_{\bA_2^-}a(\infty)=\pi_{\bA_2^-}\pi_\OO a_0=\pi_{\bA_2^-}\pi_{\bN_2^-}\pi_\OO a_0=\pi_{\bA_2^-}\pi_{\bN_1^-}\pi_\OO a_0.$$
Since $a_0\in \bP_1$, we have $\pi_{\bN_1^-}a_0=\pi_{\bA_1^-}a_0$. Thus, plugging the above expression for
$\pi_{\bA_2^-}a(\infty)$ into \eqref{pi-A1-A2-si-eq} we get
$$\pi_{\bA_1^-}a_0=\si^{-1}\pi_{\bA_2^-}\pi_\OO\pi_{\bA_1^-}a_0,$$
or equivalently,
$$(1-\th^-)\pi_{\bA_1^-}a_0=0.$$
Since $\th^-$ is nilpotent by Lemma \ref{nilpotent-operator-lem}, we conclude that $\pi_{\bA_1^-}a_0=0$
and hence $\pi_{\bA_2^-}a(\infty)=0$.

Next, we similarly analyze the relation
\begin{equation}\label{pi-A2-A1-si-eq}
\pi_{\bA_2^+}a(\infty)=\si\pi_{\bA_1^+}a_0.
\end{equation}
Namely, using \eqref{B1+a0-relation} we get
$$\pi_{\bA_1^+}a_0=\pi_{\bA_1^+}\pi_{\OO}\pi_{\bA_2^+}a(\infty),$$
so \eqref{pi-A2-A1-si-eq} can be rewritten as
$$(1-\th^+)\pi_{\bA_2^+}a(\infty)=0,$$
Hence, by nilpotency of $\th^+$, we deduce that $\pi_{\bA_2^+}a(\infty)=0$. 
Since $a(\infty)\in \bB_2^- + \bA_2^+$, we obtain that
$$\pi_{\bN_2^+}a(\infty)=\pi_{\bA_2^+}a(\infty)=0.$$

Now using \eqref{B1+a0-relation} together with the fact
that $a_0\in \bB_1^+ + \bA_1^-$ we deduce that 
$$a_0=\pi_{\bB_1^+}a_0=\pi_\bH a(\infty).$$
In particular, $a_0\in\bH$. Since $a(\infty)=a_0+a_1$, we get
that 
$$\pi_{\bN_2^+}a_1=\pi_{\bN_2^+}a(\infty)=0.$$ 
But $a_1\in\bM(\OO(1))\sub \bN_2^+$, so we deduce that $a_1=0$.

Thus, $a(\infty)=a_0\in\bH$ satisfies $a_0=\si^{-1}a_0$. Since $\si$ is transitive, this is possible only
when $a_0$ is a scalar matrix. This proves that $h^0(\AA)=1$, so we can apply Proposition \ref{parabolic-constr-prop}.
\ed

\begin{rem}
The order $\AA(\VV,I(\bP_1,\bP_2,\sigma))$ constructed above is an example of a {\it noncommutative nodal order}, the notion defined
and studied in \cite{BurD}.
\end{rem}

\subsection{Computation of the solution of the AYBE}\label{AYBE-comp-sec}

Now we are ready to prove our first main result. Note that although the classification of trigonometric solutions was obtained over $\C$,
below we work over any algebraically closed field of characteristic $0$.

\begin{thm}\label{trig-AYBE-thm} Let $\AA$ be the spherical order on $C$ constructed in Sec.\ \ref{BD-sph-order-sec}
from associative BD data $(\si,\Ga_1,\Ga_2)$.
Then the solution of the AYBE associated with the pair of families 
$$(\AA\ot \LL), (V\ot\OO_p),$$ 
where $\LL$ runs through degree $0$ line bundles on $C$,
is equivalent to
the trigonometric solution \eqref{r-gen-AYBE-sol} associated with the data $(\si^{-1},\Ga_2,\Ga_1)$.
\end{thm}

\Pf . {\bf Overview of the calculation}.
Recall that line bundles of degree $0$ on $C$ are parametrized by $\G_m$:
for each $\la\in k^*$ we have a line bundle $L_\la$ obtained by descending the trivial line bundle on $\P^1$
and using the identification of fibers 
$$\OO_{\P^1}|_0=k\rTo{\cdot \la} k=\OO_{\P^1}|_\infty.$$
In particular, we have natural trivializations of $L_\la$ at smooth points of $C$, and natural isomorphisms
$$L_\la^{-1}\ot L_{\la'}\simeq L_{\la^{-1}\la'}.$$
Applying Proposition \ref{AYBE-order-sol-prop} in our setting, we see
that the for generic $\la_1,\la_2\in k^*$ and generic points $x,y$ on $C$,
we have
\begin{equation}
r^{\la_1,\la_2}_{x,y}=\sum_{i,j} \phi_{x,y}^{\la_1^{-1}\la_2}(e_{ij})\ot e_{ji},
\end{equation}
where 
$$\phi_{x,y}^{\la}:=\ev_y\circ r_x^{-1}:\Mat_n(k)\to \Mat_n(k),$$
with the maps $r_x$ and $\ev_y$ given by
$$r_x:H^0(\AA\ot L_\la(x))\to \Mat_n(k): a\mapsto \res_x(\frac{dt}{t}\cdot a),$$
$$\ev_y:H^0(\AA\ot L_\la(x))\to \Mat_n(k): a\mapsto a(y),$$
where we use the trivializations of $L_\la$ at $x$ and $y$.

Thus, our goal is to calculate the map $\phi=\phi_{x,y}^\la$. To this end we will fix
$a\in H^0(\AA\ot L_\la(x))$ and set $b:=r_x(a)$. In the first part of the proof we will
express $a$ in terms of $b$, and then evaluate the result at $y$. This gives some expression for $\phi$, and hence
for the $r$-matrix. In the remaining part of the proof we will show
that the obtained expression is equivalent to the formula \eqref{r-gen-AYBE-sol}.

\medskip

\noindent
{\bf Step 1. The components of $a$}.
For $x\in\P^1$, $x\neq 0,\infty$, the space $H^0(\AA\ot L_\la(x))$ can be identified with the space 
of $n\times n$ matrices $a=(a_{ij})$
with $a_{ij}\in H^0(\OO(\eps_i-\eps_j)(x))$, such that 
\begin{equation}\label{a-isotropic-condition}
(\la\cdot a(0),a(\infty))\in I(\bP_1,\bP_2,\si).
\end{equation}

We identify rational sections of $\OO(1)$ (resp., $\OO$, resp., $\OO(-1)$) with rational functions of $t$ (the coordinate
on $\P^1\setminus\{\infty\}$) of degree $\le 1$ (resp., $\le 0$, resp., $\le -1$).
Since $a$ can have a pole of order at most $1$ at the point $x\in \P^1$, we can write
$$a=\frac{-xa_0+a_1t+a_2t^2}{t-x},$$
where 
$$\pi_{\OO(-1)}a_1=0, \ \ a_2\in\bM(\OO(1)).$$
This implies that
\begin{equation}\label{b-res-a-eq}
b:=\Res_{t=x} a\cdot \frac{dt}{t}=-a_0+a_1+a_2x,
\end{equation}
while
$$a(0)=a_0\in \bP_1, \ \ a(\infty)=-x\pi_{\OO(-1)}a_0+\pi_\OO a_1+a_2\in \bP_2.$$
Due to the description \eqref{Levi-P-eq} of the Levi components, the condition 
\eqref{a-isotropic-condition} is equivalent to 
\begin{equation}\label{a-isotropic-H-condition}
\la\cdot\pi_\bH a(0)=\si^{-1}\pi_\bH a(\infty),
\end{equation}  
\begin{equation}\label{a-isotropic-pm-conditions}
\la\cdot\pi_{\bA_1^\pm} a(0)=\si^{-1}\pi_{\bA_2^\pm} a(\infty).
\end{equation}

\medskip

\noindent
{\bf Step 2. Analyzing $a(\infty)$}.
We look at the formula $a(\infty)=-x\pi_{\OO(-1)}a_0+\pi_\OO a_1+a_2$ and combine the conditions
$a(\infty)\in \bP_2$ with $a_2\in \bM(\OO(1))$.
We get $\pi_\OO a(\infty)=\pi_\OO a_1$, and since $\pi_{\bB_1^+}\pi_{\OO(1)}=0$,
we obtain
\begin{equation}\label{pi-B1+a1-eq}
\pi_{\bB_1^+}a_1=\pi_{\bB_1^+}\pi_{\OO}a_1=
\pi_{\bB_1^+}\pi_\OO a(\infty)=
\pi_{\bB_1^+}\pi_\OO \pi_{\bB_2^+}a(\infty)=
\pi_{\bB_1^+}\pi_\OO \pi_{\bA_2^+}a(\infty)+\pi_\bH a(\infty)
\end{equation}
Since $a_2\in\bM(\OO(1))$, we deduce that
\begin{equation}\label{pi-B1+a0-eq}
\pi_{\bB_1^+}a_0=\pi_{\bB_1^+}(-b+a_1)=
\pi_{\bB_1^+}(-b+\pi_\OO\pi_{\bA_2^+}a(\infty))+\pi_{\bH} a(\infty).
\end{equation}

\medskip

\noindent
{\bf Step 3. Analyzing one of the relations \eqref{a-isotropic-pm-conditions}}.
Next, we analyze the relation
$$\pi_{\bA_1^-} a_0=\la^{-1}\si^{-1}\pi_{\bA_2^-} a(\infty).$$
Since $\pi_{\bA_2^-}\pi_{\OO(1)}=0$, we have
$$\pi_{\bA_2^-} a(\infty)=\pi_{\bA_2^-}(-x\pi_{\OO(-1)}a_0+\pi_\OO a_1).$$
Note that \eqref{b-res-a-eq} implies 
$$\pi_{\OO(-1)}a_0=-\pi_{\OO(-1)}b, \ \ -\pi_\OO a_0+\pi_\OO a_1=\pi_\OO b.$$
Thus, we have
$$\pi_{\bA_2^-} a(\infty)=\pi_{\bA_2^-}(x\pi_{\OO(-1)}b+\pi_\OO b+\pi_\OO a_0).$$
Further, we have $a_0\in\bB_1^+ +\bA_1^-$, so
$$\pi_{\bA_2^-}\pi_\OO a_0=\pi_{\bA_2^-}\pi_{\bN_2^-}\pi_\OO a_0=\pi_{\bA_2^-}\pi_{\bN_1^-}\pi_\OO a_0=
\pi_{\bA_2^-}\pi_{\bA_1^-}\pi_\OO a_0.$$
Thus, we obtain
$$\pi_{\bA_1^-}a_0=\la^{-1}\si^{-1}\pi_{\bA_2^-}(x\pi_{\OO(-1)}b+\pi_\OO b+\pi_\OO\pi_{\bA_1^-}a_0),$$
or equivalently,
\begin{equation}\label{pi-A1-a0-eq}
\pi_{\bA_1^-}a_0=(1-\la^{-1}\th^-)^{-1}\la^{-1}\si^{-1}\pi_{\bA_2^-}(x\pi_{\OO(-1)}b+\pi_\OO b).
\end{equation}

\medskip

\noindent
{\bf Step 4. Analyzing the other relation in \eqref{a-isotropic-pm-conditions}}.
Similarly, we analyze the relation
$$\pi_{\bA_2^+}a(\infty)=\la\cdot\si\pi_{\bA_1^+}a_0.$$
Using \eqref{pi-B1+a0-eq} we can rewrite this as
$$\pi_{\bA_2^+}a(\infty)=\la\cdot\si\pi_{\bA_1^+}(-b+\pi_\OO\pi_{\bA_2^+}a(\infty))$$
and we get
$$\pi_{\bA_2^+}a(\infty)=-\la(1-\la\th^+)^{-1}\si\pi_{\bA_1^+}b.$$
Using \eqref{pi-B1+a0-eq} again, we get
\begin{equation}\label{pi-N1+a0-eq}
\pi_{\bN_1^+}a_0=-\pi_{\bN_1^+}b-\la\pi_\OO(1-\la\th^+)^{-1}\si\pi_{\bA_1^+}b.
\end{equation}

\medskip

\noindent
{\bf Step 5. Consequence of \eqref{a-isotropic-H-condition} and expressing $a_0$ in terms of $b$}. 
Let us combine the equation $\pi_\bH a(\infty)=\la\si\pi_\bH a_0$
with the relation
$$\pi_\bH a_0=-\pi_\bH b+\pi_\bH a(\infty)$$
that follows by applying $\pi_\bH$ to \eqref{pi-B1+a0-eq}. We get
$$\pi_\bH a_0=-\pi_\bH b+\la\si\pi_\bH a_0,$$
so that
\begin{equation}\label{pi-Ha0-eq}
\pi_\bH a_0=-(1-\la\si)^{-1}\pi_\bH b
\end{equation}
(we assume that $\la$ is generic, so $1-\la\si$ is invertible).
Note that since $a_0\in \bB_1^+ +\bA_1^-$, formulas \eqref{pi-A1-a0-eq}, \eqref{pi-N1+a0-eq}
and \eqref{pi-Ha0-eq} completely determine $a_0$ in terms of $b$.

\medskip

\noindent
{\bf Step 6. Expressing $a_2$ and $a_1$ in terms of $b$}. 
We have 
$$a_2=\pi_{\OO(1)}a(\infty)=\pi_{\OO(1)}\pi_{\bA_2^+}a(\infty),$$
so using \eqref{pi-N1+a0-eq} we get
\begin{equation}\label{a2-formula}
a_2=-\la\pi_{\OO(1)}(1-\la\th^+)^{-1}\si\pi_{\bA_1^+}b.
\end{equation}
Thus, we expressed $a_0$ and $a_2$ in terms of $b$.
On the other hand, we have
$$a_1=b+a_0-xa_2,$$
which allows to express $a_1$ in terms of $b$ as well.

\medskip

\noindent
{\bf Step 7. The formula for $\phi$}.
To compute $\phi$, we have to evaluate $a$ at another point $y$ and express
the result in terms of $b$:
$$a(y)=-\frac{x}{y-x} a_0+\frac{y}{y-x} a_1+ \frac{y^2}{y-x} a_2=
a_0+\frac{y}{y-x} b +ya_2.$$

Taking into account \eqref{pi-A1-a0-eq}, \eqref{pi-N1+a0-eq}, \eqref{pi-Ha0-eq} and \eqref{a2-formula}, we obtain
$$a(y)=\phi(b):=\phi_0(b)+\phi_-(b)-\phi_+(b)+x\cdot\psi_-(b)-y\cdot\psi_+(b),$$
where
$$\phi_0(b)=\frac{y}{y-x} b-\pi_{\bN_1^+}b -(1-\la\si)^{-1}\pi_\bH b,$$
$$\phi_-(b)=(1-\la^{-1}\th^-)^{-1}\la^{-1}\si^{-1}\pi_{\bA_2^-}\pi_\OO b,$$
$$\phi_+(b)=\la\pi_\OO(1-\la\th^+)^{-1}\si\pi_{\bA_1^+}b,$$
$$\psi_-(b)=(1-\la^{-1}\th^-)^{-1}\la^{-1}\si^{-1}\pi_{\bA_2^-}\pi_{\OO(-1)}b,$$
$$\psi_+(b)=\la\pi_{\OO(1)}(1-\la\th^+)^{-1}\si\pi_{\bA_1^+}b.$$

In the rest of the proof we will simplify this formula and identify the corresponding associative $r$-matrix
$$r=\sum_{i,j}\phi(e_{ij})\ot e_{ji}$$
with \eqref{r-gen-AYBE-sol}.

\medskip

\noindent
{\bf Step 8. Expanding the expressions with $\th^\pm$}.
First, we claim that
$$\pi_{\bA_1^+}\pi_{\OO}\si\pi_{\bA_1^+}=\pi_{\bA_1^+}\si\pi_{\bA_1^+}.$$
Indeed, since $\si(\a)\in \Pi_2$, this follows immediately from the inclusion 
$$\Pi_1\cap \Pi_2\sub M(\OO).$$
Hence, we have
$$(1-\la\th^+)^{-1}\la\si\pi_{\bA_1^+}=\sum_{m\ge 1}\la^m(\si\pi_{\bA_1^+})^m.$$

Similarly, we see that
$$\pi_{\bA_2^-}\pi_{\OO}\si^{-1}\pi_{\bA_2^-}=\pi_{\bA_2^-}\si^{-1}\pi_{\bA_2^-},$$
hence,
$$(1-\la^{-1}\th^-)^{-1}\la^{-1}\si^{-1}\pi_{\bA_2^-}=\sum_{m\ge 1}\la^{-m}(\si^{-1}\pi_{\bA_2^-})^m.$$

\medskip

\noindent
{\bf Step 9. Matching the part with $\phi_0$ with $r_{\si^{-1},\const}$}.
Let us set $z=x/y$. Then we have
$$\phi_0(b)=\frac{1}{1-z} b-\pi_{\bN_1^+}b -(1-\la\si)^{-1}\pi_\bH b.$$
Hence, 
$$\phi_0(e_\a)=\begin{cases} \frac{z}{1-z} e_\a, & \a>0, \\ \frac{1}{1-z} e_\a, & \a<0, \end{cases}$$
while 
$$\phi_0(h_i)=\frac{1}{1-z} h_i-\sum_{m\ge 0}\la^mh_{\si^m(i)}=\frac{1}{1-z} h_i-
(1-\la^n)^{-1}\sum_{k=0}^{n-1}\la^k h_{\si^k(i)}.$$
Hence,
\begin{align*}
&\sum_{i,j}\phi_0(e_{ij})\ot e_{ji}=\frac{z}{1-z}\sum_{\a>0} e_{\a}\ot e_{-\a}+\frac{1}{1-z}\sum_{\a>0} e_{-\a}\ot e_{\a}\\
&+\frac{1}{1-z}\sum_i h_i\ot h_i-(1-\la^n)^{-1}\sum_{i=1}^n\sum_{k=0}^{n-1} \la^k h_{\si^k(i)}\otimes h_i 
=r_{\const}(\si^{-1},\la,z).
\end{align*}

\medskip

\noindent
{\bf Step 10. Matching the parts depending on $\Pi_1$}.
Next, we have $\pi_{A_2^-}\pi_{\OO}e_{-\a}\neq 0$ precisely when $\a\in \Pi_2$ and $\a>0$. Thus,
$$\sum_{i,j}\phi_-(e_{ij})\ot e_{ji}=
\sum_{\a\in\Pi_2, \a>0, k\ge 1}\la^{-k}(\si^{-1}\pi_{\bA_2^-})^ke_{-\a}\ot e_{\a}.$$

Recall that the formula \eqref{r-gen-AYBE-sol} for the data $(\si^{-1},\Ga_2,\Ga_1)$
uses the partially defined operation $\tau^{-1}$ which is defined only on $\Pi_2$ and is given by $\tau^{-1}(\a)=\si^{-1}(\a)$.
Thus, we have
$$(\si^{-1}\pi_{\bA_2^-})^ke_{-\a}=\begin{cases} e_{-\tau^{-k}(\a)}, & \tau^k(\a) \text{ is defined},\\ 0, & \text{otherwise}.\end{cases}$$
Thus, we can rewrite the above formula as
$$\sum_{i,j}\phi_-(e_{ij})\ot e_{ji}=
\sum_{\a\in\Pi_2, \a>0, k\ge 1}\la^{-k}e_{-\tau^{-k}(\a)}\ot e_{\a}.$$

Similarly, $\pi_{A_2^-}\pi_{\OO(-1)}e_{-\a}\neq 0$ precisely when $\a\in \Pi_2$ and $\a<0$.
Thus,
$$\sum_{i,j}\psi_-(e_{ij})\ot e_{ji}=
\sum_{\a\in\Pi_2, \a<0, k\ge 1}\la^{-k}(\si^{-1}\pi_{\bA_2^-})^ke_{-\a}\ot e_{\a}=
\sum_{\a\in\Pi_2, \a<0, k\ge 1}\la^{-k}e_{-\tau^{-k}(\a)}\ot e_{\a}.$$

Next, again using the fact that for $\b\in\Pi_2$ one has $\b\in M(\OO)$ if and only if $\b>0$, we see that
for $k\ge 1$ the condition $\pi_{\OO}(\si\pi_{A_1^+})^ke_{\a}\neq 0$ is equivalent to
$$\a\in\Pi_1, \si(\a)\in\Pi_1, \ldots, \si^{k-1}(\a)\in\Pi_1, \text{ and } \si^k(\a)>0.$$
Equivalently, $\b=\si^k(\a)\in\Pi_2$ satisfies, $\b>0$ and $\tau^{-k}(\b)$ is defined.
Hence,
\begin{align*}
&\sum_{i,j}\phi_+(e_{ij})\ot e_{ji}=
\sum_{\a\in\Pi_1, k\ge 1, \si^k(\a)>0}\la^k(\si\pi_{\bA_1^+})^ke_{\a}\ot e_{-\a}=\\
&\sum_{\b\in\Pi_2,\b>0, k\ge 1}\la^k e_\b\ot (\si^{-1}\pi_{\bA_2^-})^ke_{-\b}=
\sum_{\b\in\Pi_2,\b>0, k\ge 1}\la^k e_\b\ot e_{-\tau^{-k}(\b)}.
\end{align*}

Similarly, for $k\ge 1$, the condition $\pi_{\OO(1)}(\si\pi_{A_1^+})^ke_{\a}\neq 0$ is equivalent to 
$$\a\in\Pi_1, \si(\a)\in\Pi_1, \ldots, \si^{k-1}(\a)\in\Pi_1, \text{ and } \si^k(\a)<0.$$
Thus, we can rewrite
\begin{align*}
&\sum_{i,j}\psi_+(e_{ij})\ot e_{ji}=
\sum_{\a\in\Pi_1, k\ge 1, \si^k(\a)<0}\la^k(\si\pi_{\bA_1^+})^ke_{\a}\ot e_{-\a}\\
&=\sum_{\b\in\Pi_2,\b<0, k\ge 1}\la^k  e_\b\ot (\si^{-1}\pi_{\bA_2^-})^ke_{-\b}=
\sum_{\b\in\Pi_2,\b<0, k\ge 1}\la^k  e_\b\ot e_{-\tau^{-k}(\b)}.
\end{align*}

Combining the above computations we obtain
$$r(\la;x,y)=r_{\si^{-1},\Ga_2,\Ga_1}(\la;x,y).$$
\ed

\begin{rem}
1. As shown in \cite{LP-YB}, the $A_\infty$-category split generated by a pair of $1$-spherical objects is completely determined by
the corresponding formal solution of the AYBE. Since the pair of objects $(\AA,V\ot\OO_x)$ split generates the perfect derived category
of right $\AA$-modules $\Perf(\AA)$ (see \cite[Lem.\ 3.2.1]{P-ainf-orders}), this means that whenever we have a pair of $1$-spherical
objects in some $A_\infty$-category $\CC$, giving rise to a nondegenerate trigonometric solution of the AYBE, we get a fully faithful 
embedding of $\Perf(\AA)$ into $\CC$. In particular, this applies to a pair of the form $(\VV,\OO_x)$ on the wheel of $n$ projective lines $G_n$
(aka {\it standard $n$-gon}), where $\VV$ is a simple vector bundle on $G_n$, since as was shown in \cite{P-nodal-massey}
such a pair gives rise to the trigonometric solution corresponding to some BD data $(\si_0,\si,\Ga_1,\Ga_2)$ with commuting $\si_0$ and $\si$.
This suggests that the corresponding spherical order $\AA$ should be isomorphic to $p_*\End(\VV)$, where $p:G_n\to C$
is the natural morphism contracting all the components not containing the point $x\in G_n$. 

\noindent
2. For each associative BD data we constructed in \cite{LP-YB} a pair of $1$-spherical objects in the Fukaya category $\FF(\Sigma)$ of a certian square-tiled 
(noncompact) surface $\Sigma$, giving rise to the corresponding trigonometric solution of the AYBE. 
It follows that there is a fully faithful functor from $\Perf(\AA)$ to $\FF(\Sigma)$. It is plausible that it is in fact an equivalence. In the case when $\Sigma$
has genus $1$, this follows from \cite[Prop.\ 2.3.6]{LP-YB}.
\end{rem}

\subsection{More on symmetric spherical orders}\label{more-sph-sec}

Recall that an order $\AA$ over an integral projective curve $C$ is called {\it weakly spherical} if $h^0(C,\AA)=h^1(C,\AA)=1$.

Part (i) of the following result is somewhat analogous to \cite[Lem.\ 3.16]{BG}.

\begin{prop}\label{w-spherical-order-prop}
Let $\AA$ be a weakly spherical order of rank $n^2$ over an integral projective curve $C$, where $n\ge 2$,
and let $\nu: \wt{C}\to C$ be the normalization map.

\noindent
(i) The genus of $\wt{C}$ is $\le 1$. Furthermore, if $\wt{C}$ is of genus $1$ then there exists
a simple vector bundle $\VV$ on $\wt{C}$ such that $\AA\simeq \nu_*\und{\End}(\VV)$.

\noindent
(ii) Assume that $\wt{C}\simeq\P^1$. Then there exists a vector bundle $\VV$ on $\P^1$ of one of the two types,
\begin{enumerate}
\item
$\VV=\OO^{\oplus m}\oplus \OO(1)^{\oplus (n-m)}$,
\item
$\VV=\OO\oplus \OO(1)^{\oplus (n-2)}\oplus \OO(2)$,
\end{enumerate}
and an embedding of orders $\AA\sub \nu_*\und{\End}(\VV)$ such that the quotient is a torsion sheaf $T$ of length
$n^2$. 
Conversely, any such suborder $\AA$ with $h^0(\AA)=1$ is weakly spherical.
\end{prop}

\Pf . (i) Set $\wt{\AA}:=\nu^*\AA/\tors$. Then $\AA$ is a subsheaf of $\nu_*\wt{\AA}$ and the quotient $\nu_*\wt{\AA}/\AA$ is a torsion sheaf.
This implies that $h^1(\wt{\AA})=h^1(\nu_*\wt{\AA})\le 1$.
Let $\BB$ be a maximal order on $\wt{C}$ such that $\wt{\AA}\sub \BB$ and the quotient is a torsion sheaf.
Such a maximal order always exists (see e.g., \cite[Prop.\ 4.5]{Chan}). Since $k$ is algebraically closed, 
there exists a vector bundle $\VV$ on $\wt{C}$ such that
$\BB\simeq \und{\End}(\VV)$ (see \cite[Cor.\ 13.2]{Chan}). Note that $h^1(\BB)\le h^1(\wt{\AA})\le 1$ since $\BB/\wt{\AA}$ is torsion.
Thus, $\dim\Ext^1(\VV,\VV)\le 1$. By Serre duality, this is equivalent to $\dim \Hom(\VV,\VV\ot\om_{\wt{C}})\le 1$.
This is possible only when $\wt{C}$ is rational or an elliptic curve.

Assume now that $\wt{C}$ is an elliptic curve. Then we get that $\VV$ is a simple vector bundle, so $h^0(\BB)=h^1(\BB)=1$.
Now we have an exact sequence
\begin{equation}\label{orders-ex-seq}
0\to \AA\to \nu_*\und{\End}(\VV)\to T\to 0
\end{equation}
where $T$ is a torsion sheaf. Then
$$\ell(T)=\chi(T)=\chi(\und{\End}(\VV))-\chi(\AA)=0,$$
so $T=0$.


\noindent
(ii) As in part (i), we have an embedding $\wt{\AA}\sub\BB=\und{\End}(\VV)$, where $h^1(\BB)\le 1$. The condition 
$\dim \Hom(\VV,\VV(-2))\le 1$ immediately implies that up to tensoring with a line bundle, $\VV$ is of one of the types
(1) or (2). 

We again have an exact sequence \eqref{orders-ex-seq}, which gives
$$\ell(T)=\chi(T)=\chi(\und{\End}(\VV))-\chi(\AA)=\chi(\und{\End}{\VV})=n^2.$$

Conversely, assume we have a suborder $\AA\sub \nu_*\und{\End}(\VV)$ such that $h^0(\AA)=1$ and the quotient has length $n^2$.
Then we get $\chi(\AA)=0$, so $h^1(\AA)=1$.
\ed

Now let us specialize to the case when $C$ is the irreducible nodal curve of arithmetic genus $1$ and the case
of symmetric spherical orders.

\begin{prop} Let $\AA$ be a symmetric spherical order of rank $n^2$ on the irreducible nodal curve $C$ of arithmetic genus $1$.
Then
$\AA$ arises by the construction of Lemma \ref{sym-sph-order-isotr-lem} from some vector bundle $\VV$ of rank
$n$ over $\P^1$ and a maximal isotropic subalgebra $I\sub \End(\VV|_0)\oplus \End(\VV|_\infty)$.
Furthermore, replacing $\VV$ by $\VV(i)$, we can achieve that $\VV\simeq \OO^{\oplus m}\oplus \OO(1)^{\oplus (n-m)}$.
\end{prop}

\Pf . First, as in Proposition \ref{w-spherical-order-prop}(ii), we get an embedding of orders
$\AA\sub \nu_*\und{\End}(\VV)$, with the quotient $T$ of length $n^2$.
Next, we are going to use the nonzero map $\tau:\AA\to \OO_C$, such that $\tau(xy)=\tau(yx)$.
We claim that after rescaling $\tau$ by a nonzero constant, we have a commutative diagram
\begin{equation}\label{tau-tr-diagram}
\begin{diagram}
\AA &\rTo{} &\nu_*\und{\End}(\VV)\\
\dTo{\tau}&&\dTo{\tr}\\
\OO_C&\rTo{} & \nu_*\OO_{\P^1}
\end{diagram}
\end{equation}
Indeed, looking at the generic point we see that $\tr=f\cdot\tau$ for some nonzero rational function $f$ on $C$.
Comparing the values on $1\in H^0(\AA)$ we get
$$n=f\cdot \tau(1).$$
But $\tau(1)$ is a global function on $\OO_C$, hence a constant. Therefore, $f$ is also a constant.

Now let $q\in C$ denote the singular point, so that $\nu^{-1}(q)=\{0,\infty\}$, and let
$$I\sub \nu_*\und{\End}(\VV)|_q=\End(\VV|_0)\oplus \End(\VV|_\infty)$$
denote the image of $\AA|_q$. Clearly this is a subalgebra.
By the commutativity of diagram \eqref{tau-tr-diagram}, we get
that $\tr(a)=\tr(b)$ for $(a,b)\in I$. Hence, $I$ is an isotropic subalgebra,
in particular, $\dim I\le n^2$.
Hence, the exact sequence
$$0\to I\to \End(\VV|_0)\oplus \End(\VV|_\infty)\to T/\mg_qT\to 0$$
shows that $\dim T/\mg_qT\ge n^2$. But $T$ has length $\le n^2$, so this is possible only
if $\mg_qT=0$, $\ell(T)=n^2$ and $\dim I=n^2$. Thus, $I$ is maximal isotropic and $\AA=\AA(\VV,I)$, and
twisting $\VV$, we can achieve that
$\VV\simeq \OO^{\oplus m}\oplus \OO(1)^{\oplus (n-m)}$.
\ed

\section{Sheaves of Lie algebras and the CYBE}

\subsection{Formal solutions of the CYBE and Manin triples}\label{formal-CYBE-sec}

The results of this subsection are well known to the experts. We refer to \cite[ch.\ 6,7]{ES} for some background.
For a Lie algebra $\fg$ over $\C$, we consider the classical Yang-Baxter equation of the form
\begin{equation}\label{CYBE-eq}
[r^{12}(x_1,x_2),r^{13}(x_1,x_3)]+[r^{12}(x_1,x_2),r^{23}(x_2,x_3)]+[r^{13}(x_1,x_3),r^{23}(x_2,x_3)]=0,
\end{equation}
where $r(x_1,x_2)$ takes values in $\fg\ot \fg$. Classically $r(x_1,x_2)$ is viewed as a meromorphic function where $x_1$ and $x_2$
vary in some open domain in $\C$ (see \cite{BD-genYB}). In the case when $r$ depends only on the difference $x_1-x_2$,
we get the equation \eqref{CYBE-intro-eq}.
If $\fg$ is a Lie algebra over a field $k$ then one can also assume that $x_1$ and $x_2$
vary in a smooth curve $U$ over $k$, and assume that $r(x_1,x_2)$ is a rational function on $U\times U$.
The CYBE is often coupled with the unitarity condition
\begin{equation}\label{CYBE-unitary-eq}
r^{21}(x_2,x_1)=-r(x_1,x_2).
\end{equation}

For a Lie algebra $\fg$ over a field $k$ it makes sense to consider the following formal solutions of \eqref{CYBE-eq} of the form
$$r(x_1,x_2)\in \fg\ot \fg(\!(x_1)\!)[\![x_2]\!]$$ 
and impose the equation \eqref{CYBE-eq} in 
$$U(\fg)^{\ot 3}\ot R,$$
where $R=k(\!(x_1)\!)(\!(x_2)\!)[\![x_3]\!]$.
Indeed, to make sense of the equation we use the natural embeddings
$$k(\!(x_2)\!)\sub k(\!(x_1)\!)(\!(x_2)\!), \ \ k(\!(x_1)\!)\sub k(\!(x_1))(\!(x_2)\!), \ \ k(\!(x_1)\!)[\![x_2]\!]\sub k(\!(x_1)\!)(\!(x_2)\!)$$
and view each $r^{ij}(x_i,x_j)$, for $i<j$ as an element of $U(\fg)^{\ot 3}\ot R$.

Let us assume in addition that $\fg$ is a finite dimensional and is equipped with an invariant nondegenerate symmetric pairing $(\cdot,\cdot)$, 
and let $\Om\in \fg\ot \fg$
be the corresponding Casimir element. 
Let us consider $r$ of the form
\begin{equation}\label{r-matrix-Om-form-eq}
r(t,u)=\frac{\Om}{t-u}+r_{\reg}, \ \ r_{\reg}\in \fg\ot\fg[\![t,u]\!],
\end{equation}
where $\frac{1}{t-u}:=\sum_{m\ge 0}t^{-m-1}u^m\in k(\!(t)\!)[\![u]\!]$.
Then we can impose the unitarity condition
$$r^{21}_{\reg}(t,u)=-r_{\reg}(u,t).$$

To a formal unitary solution $r$ of \eqref{CYBE-eq} of the form \eqref{r-matrix-Om-form-eq} we can associate a Lie subalgebra
$$\fg(r)\sub \fg(\!(t)\!)$$ 
as follows. Consider the expansion
$$r^{21}(t,u)=\sum_{n\ge 0}c_n(t)u^n\in \fg\ot \fg(\!(t)\!)[\![u]\!],$$
where $c_n\in \fg\ot \fg(\!(t)\!)$.
Then
$$\fg(r):=\operatorname{span}\{(\phi\ot \id)c_n(t) \ |\ \phi\in \fg^*,n\ge 0\}\sub \fg(\!(t)\!).$$
It is clear that $\fg(r)$ is a complementary subspace to $\fg[\![t]\!]$. The fact that it is a Lie subalgebra follows from the CYBE
(see Proposition \ref{formal-prop} below).

Let us equip $\fg(\!(t)\!)$ with the invariant nondegenerate form
$$\lan X(t),Y(t)\ran:=\Res_0((X(t),Y(t))dt).$$


\begin{prop}\label{formal-prop}
In the above situation $\fg(r)$ is an isotropic Lie subalgebra of $\fg(\!(t)\!)$ and 
$$\fg(\!(t)\!)=\fg(r)\oplus \fg[\![t]\!].$$
In other words, $(\fg(\!(t)\!),\fg(r),\fg[\![t]\!])$ is an (infinite-dimensional) Manin triple.
Furthermore, the above construction establishes a bijection between the set of formal unitary solutions of the CYBE of the form \eqref{r-matrix-Om-form-eq}
and Manin triples of the form $(\fg(\!(t)\!),L,\fg[\![t]\!])$.
\end{prop}

\Pf . Using the same arguments as in \cite[Prop.\ 6.2]{ES}, one checks that starting from a formal solution $r$ of the CYBE 
of the form \eqref{r-matrix-Om-form-eq} we get a Lie subalgebra $\fg(r)\sub \fg(\!(t)\!)$. 
Furthermore, the condition that $\fg(r)$ is isotropic is equivalent to $r$ being unitary.

Conversely, given a Lie subalgebra $L\sub \fg(\!(t)\!)$, complementary to $\fg[\![t]\!]$, we define 
$c_n\in \fg\ot \fg(\!(t)\!)$, for $n\ge 0$, from the condition
$$c_n=\sum X_i\ot c_n^i,$$
where $(X_i)$ is an orthonormal basis of $\fg$, and $c_n^i$ is the unique element of $L$, such that
$$c_n^i\equiv \frac{X_i}{t^{n+1}} \mod \fg[\![t]\!].$$
Then for $r(t,u)=\sum_{n\ge 0}c^{21}_n(t)u^n$ the CYBE is equivalent to the vanishing of 
\begin{align*}
&\CYBE(r)=\sum_{i,j,n,m}[c_n^i(x_1),c_m^j(x_2)]\ot X_i\ot X_j \cdot x_2^nx_3^m+\\
&\sum_{i,j,n,m}c_n^i(x_1)\ot [X_i,c^j_m(x_2)]\ot X_j \cdot x_2^nx_3^m+
\sum_{i,j,n,m}c_n^i(x_1)\ot c_m^j(x_2)\ot [X_i,X_j]\cdot x_3^{m+n}.
\end{align*}
It is easy to check that in fact
$\CYBE(r)$ is regular, i.e., belongs to $\fg^{\ot 3}[\![x_1,x_2,x_3]\!]$.
On the other hand, the above formula shows that
$$\CYBE(r)\in L\ot \fg^{\ot 2}[\![x_2,x_3]\!],$$
where we view $L$ as a subspace of $\fg(\!(x_1)\!)$.
Since $L\cap \fg[\![x_1]\!]=0$, this is possible only if $\CYBE(r)=0$.

It is easy to see that the above two constructions are inverses of each other.
\ed

\subsection{Spherical sheaves of Lie algebras and Manin triples}\label{spherical-Lie-sec}

Let $C$ be an integral projective curve over $k$.
We consider the following Lie analog of Definition \ref{sph-order-def}.

\begin{defi}\label{sph-Lie-def}
Let $\LL$ be a coherent sheaf of Lie algebras on $C$. We say that $\LL$ is {\it symmetric spherical} if 
$\LL$ is acyclic, i.e., $H^*(C,\LL)=0$, and there
exists a symmetric $\OO$-bilinear invariant pairing
$$\kappa:\LL\ot \LL\to \om_C$$
inducing an isomorphism
$$\nu_\kappa:\LL\to \und{\Hom}(\LL,\om_C):x\mapsto (y\mapsto \kappa(x,y)).$$
The invariance of $\kappa$ means
$$\kappa([x,y],z)=\kappa(x,[y,z]).$$
\end{defi}

Note that $\LL$ is automatically torsion free since it has no global sections. Hence, $\LL$ is locally free over the smooth
locus of $C$. 

\begin{lem}\label{sph-Lie-def-lem} 
In Definition \ref{sph-Lie-def}, instead of requiring that $\nu_\kappa$ is an isomorphism globally, it is enough to assume that $\nu_\kappa$ is generically an isomorphism.
\end{lem}

\Pf . It is well known that for a torsion free sheaf $\FF$ on $C$ with $H^0(C,\FF)=0$ one has $\und{\Ext}^{>0}(\FF,\om_C)=0$ 
(see e.g., the proof of \cite[Prop.\ 3.2.2]{P-ainf-orders}).
Hence, by Serre duality, one has
$$\chi(\und{\Hom}(\LL,\om_C))=-\chi(\LL)=0.$$
Now suppose $\nu_\kappa$ is generically an isomorphism. Then it is injective with the quotient $Q$ which is a torsion sheaf. But we should have $\chi(Q)=0$, so $Q=0$.
\ed


We have a natural construction of Manin triples from symmetric spherical Lie algebras.

\begin{prop}\label{sph-Lie-prop}
Let $(\LL,\kappa)$ be a symmetric spherical sheaf of Lie algebras on $C$, and let $p$ be a smooth point.
We denote by $\hat{\OO}$ the completion of the local ring $\OO_{C,p}$ and by $K$ its field of fractions.
Let us consider the completion $\hat{\LL}_p$ which is a Lie algebra over $\hat{\OO}$, and let us set
$$\LL_K:=\hat{\LL}_p\ot_{\hat{\OO}}K.$$
Note that $\kappa$ induces a pairing
$$\kappa_K:\LL_K\ot_K \LL_K\to \om_{C,p}\ot K.$$
Let us equip $\LL_K$ with the symmetric bilinear form $\Res_p\circ \kappa_K$.
Then
$$(\LL_K, \LL(C-p), \hat{\LL}_p)$$
is a Manin triple. 
\end{prop}

\Pf . The nongeneracy of the form $\Res_p\circ\kappa_K$ follows from the nondegeneracy of $\kappa$ near $p$.
It is clear that $\LL(C-p)$ and $\hat{\LL}_p$ are Lie subalgebras in $\LL_K$, and that $\hat{\LL}_p$ is isotropic.
The restriction of $\kappa_K$ to $\LL(C-p)$ is given by the residue at $p$ of a section of $\om_C(C-p)$, so it is zero.

On the other hand, it is well known that the cohomology of $\LL$ can be computed by the $2$-term complex
$$\LL(C-p)\oplus \hat{\LL}_p\to \LL_K.$$
Hence, the vanishing of the cohomology of $\LL$
is equivalent to the direct sum decomposition
$$\LL_K=\LL(C-p)\oplus \hat{\LL}_p.$$
\ed

\begin{cor}\label{geom-r-matrix-constr} 
In the situation of Proposition \ref{sph-Lie-prop}, assume that 
$\fg:=\LL|_p$ is a simple Lie algebra over $k$. 
Then there exists an isomorphism of $\hat{\OO}$-linear Lie algebras
$$\hat{\LL}_p\simeq \fg\ot \hat{\OO},$$
an invariant nondegenerate symmetric pairing $(\cdot,\cdot)_{\fg}$ on $\fg$,
and an $\hat{\OO}$-generator $\eta\in \hat{\om}_{C,p}=\om_{C,p}\ot \hat{\OO}$
such that 
$$\kappa_K(X\ot f,Y\ot g)=(X,Y)_{\fg}\cdot fg\cdot \eta.$$
Thus, if we choose a formal parameter $t$ at $p$, such that $\eta=dt$,
then we get a Manin triple 
$$(\fg(\!(t)\!),\LL(C-p),\fg[\![t]\!),$$
where we identify $\LL(C-p)$ with a subspace in $\LL_K\simeq \fg(\!(t)\!)$.
\end{cor}

\Pf . The existence of a trivialization follows from the fact that $\fg$ has no nontrivial formal deformations.
The $\fg$-invariant $\hat{\OO}$-linear pairing
$$\kappa_{\hat{\OO}}:\fg\ot\fg\ot \hat{\OO}\to \hat{\om}_{C,p}$$
corresponds to a $\fg$-invariant element of $(\fg\ot \fg)^*\ot_k \hat{\om}_{C,p}$,
which necessarily has form $(\cdot,\cdot)_{\fg}\ot \eta$. 
The nondegeneracy of $(\cdot,\cdot)_{\fg}$ and the fact that $\eta$ is a generator
follow from the nondegeneracy of $\kappa$ near $p$.
\ed

Now assume that $\LL$ is a symmetric spherical sheaf of Lie algebras on $C$, equipped with a trivialization
$\LL\simeq \fg\ot \OO$ over some smooth open subset $U\sub C$, where $\fg$ is a simple Lie algebra.
There is a more direct construction of a classical $r$-matrix from this data which leads to the same Manin triple as above
(this construction is discussed in detail in \cite{BG}).  
Let us define a $\fg\ot \fg$-valued rational function $r$ on $U\times U$ with the pole of order $1$ along the diagonal as follows.
Note that the restriction of $\kappa$ to $U$ has form
$$\kappa_K(X\ot f,Y\ot g)=(X,Y)_{\fg}\cdot fg\cdot \eta$$
for a nondegenerat invariant pairing $(\cdot,\cdot)_{\fg}$ on $\fg$ and an
everywhere nonvanishing $1$-form $\eta$ over $U$.
For each $x\in U$, the map
$$\Res_x: H^0(C,\LL(x))\to \fg: s\mapsto \Res_x(s\cdot \eta)$$
is an isomorphism, so for $y\in U$, $y\neq x$, we can define a map
$$\phi(x,y):\fg\to \fg$$
as the composition of $\Res_x^{-1}$ with the evaluation at $y$.
Using the nondegenerate form $(\cdot,\cdot)_{\fg}$ on $\fg$ we convert this map into a tensor $r(x,y)\in \fg\ot \fg$:
$$r(x,y)=\sum_i \phi(x,y)(e_i)\ot e_i,$$
where $(e_i)$ is an orthonormal basis of $\fg$ with respect to $(\cdot,\cdot)_{\fg}$.
Another way to state this construction is by considering the
residue map
$$H^0(C\times U,\LL\boxtimes\LL(\De_U))\rTo{\sim} H^0(\De_U,\LL|_U\ot\LL|_U\ot \om_U^{-1})\simeq \fg\ot\fg\ot H^0(U,\om_U^{-1}),$$
where $\De_U\sub U\times U\sub C\times U$ is the diagonal divisor,
and define $r$ as the preimage of $\Om\ot \eta$, where $\Om\in \fg\ot\fg$ is the Casimir element corresponding to $(\cdot,\cdot)_{\fg}$.

The $r$-matrix $r(x,y)$ satisfies the CYBE \eqref{CYBE-eq} and the unitarity \eqref{CYBE-unitary-eq}.
Let us fix a point $p\in U$, and let $t$ be the formal parameter at $p$ such that $\eta=dt$.
We can expand $r(x,y)$ near $(p,p)$ into a series in $\fg\ot \fg(\!(t)\!)[\![u]\!]$ and consider the corresponding isotropic subalgebra $\fg(r)\sub \fg(\!(t)\!)$. 
It is easy to check that
$$\fg(r)=\LL(C-p)\sub \fg(\!(t)\!).$$


Now we will give a criterion allowing to construct a symmetric spherical sheaf of Lie algebras from a Manin triple (see Proposition
\ref{triple-to-sheaf-prop} below).
As before, we assume that $p$ is a smooth point on an integral projective curve $C$, $\hat{\OO}$ is the 
completion of $\OO_{C,p}$ and $K$ is its field of fraction.

We will use the following simple fact from geometry of singular curves.

\begin{lem}\label{sing-om-lem} 
Let $\eta$ be a rational $1$-form on $C$ such that for any $f\in \OO(C-p)$ one has $\Res_p(f\eta)=0$.
Then $\eta\in \om_C(C-p)$.
\end{lem}

\Pf . Let $\nu:\wt{C}\to C$ be the normalization. It is well known that sections of $\om_C(C-p)$ can be identified with
rational $1$-forms $\xi$ on $\wt{C}$, regular at all smooth points of $C-p$ and such that for every singular point $q\in C$
and every $\phi\in \OO_{C,q}$ one has
\begin{equation}\label{sum-residues-node-eq}
\sum_{\wt{q}\in \nu^{-1}(q)} \Res_{\wt{q}}(\phi\xi)=0.
\end{equation}

First, we claim that $\eta$ is regular at any smooth point $p'$ of $C-p$. Indeed, assume $\eta$ has a pole of order $m$ at $p'$.
For any $N>0$, we can find a function $\wt{f}\in \OO(\wt{C}-p)$ with the following properties:
\begin{itemize}
\item $f$ vanishes to order $m-1$ at $p'$; 
\item $f$ vanishes to order $N$ at all other points of $\wt{C}-p$ where $\eta$ has poles;
\item $f$ vanishes to order $N$ at preimages of all singular points of $C$. 
\end{itemize}
For sufficiently large $N$, such a function is necessarily the pull-back of a regular function
$f$ on $C-p$. But for such $f$ we will have $\Res_{p'}(f\eta)\neq 0$. Since by assumption $\Res_p(f\eta)=0$, we 
get that the sum of residues of $f\eta$ is nonzero, which is a contradiction. This proves our claim that $\eta$ is regular at
all smooth points of $C-p$.

Next, we need to check condition \eqref{sum-residues-node-eq} for $\xi=\eta$ and for every singular point $q\in C$. 
Note that this condition depends only on
$\phi$ modulo some power of the maximal ideal in $\OO_{C,q}$. Now suppose for some $q$ there exists
$\phi\in \OO_{C,q}$ such that equality \eqref{sum-residues-node-eq} does not hold. We can find a regular function $f\in \OO(C-p)$
that agrees with $\phi$ modulo sufficiently high order of the maximal ideal of $\OO_{C,q}$, and at the same time belongs to
sufficiently high power of the maximal ideal in $\OO_{C,q'}$ for every other singular point $q'$ of $C$. Since $\Res_p(f\eta)=0$,
we will again obtain that the sum of all residues of $f\eta$ is nonzero, which is a contradiction.
\ed

\begin{prop}\label{triple-to-sheaf-prop} 
Let $\fg$ be a finite dimensional Lie algebra with a nondegenerate invariant symmetric pairing $(\cdot,\cdot)$, and assume that
we have a Manin triple $(\fg\ot K,L,\fg\ot\hat{\OO})$, where $\fg\ot K$ is equipped with the pairing
$$(X\ot f,Y\ot g)=(X,Y)\cdot \Res_p(fg\eta),$$
for some nonzero rational $1$-form $\eta$ on $C$, regular and nonvanishing at $p$.

\noindent
(i) Assume that $L$ is contained in $\fg\ot K_C$, where $K_C$ is the field of rational functions on $C$, and 
that $L$ is stable under the multiplication by $\OO(C-p)$.
Then $L$ comes from a symmetric spherical sheaf of Lie algebras $\LL$ on $C$, equipped with a trivialization $\hat{L}_p\simeq \fg\ot\hat{\OO}$.

\noindent
(ii) Assume that $\a\in \Aut(\fg)$ is an automorphism of finite order $h$, preserving the pairing $(\cdot,\cdot)$, 
and $K_C\sub K_{C'}$ is a cyclic extension of degree $h$ corresponding to a cyclic covering $C'\to C$, unramified at $p$. Let us fix a point $q\in C'$ over $p$ and let us consider the induced embedding $\iota_q:K_{C'}\hra K$. We fix an identification of the Galois group of $K_{C'}/K_C$ with $\Z/h$, and
let a generator $\zeta$ of the cyclic group $\Z/h$ act on $\fg\ot K_{C'}$ by $\a\ot \zeta$.
We have a $\Z/h$-equivariant embedding into the group ring over $K$ (where $\Z/h$ acts trivially on $K$),  
$$K_{C'}\to K[\Z/h]: f\mapsto \sum_m \iota_q(\zeta^mf)\ot \zeta^{-m}$$
and the induced embedding of Lie algebras
\begin{equation}\label{twisted-Lie-alg-K-emb}
(\fg\ot K_{C'})^{\Z/h}\sub (\fg\ot K[\Z/h])^{\Z/h}\simeq \fg\ot K,
\end{equation}
(where $\zeta\in \Z/h$ acts by $\a\ot \zeta$ on $\fg\ot K[\Z/h]$).
Assume that $L$ is contained in $(\fg\ot K_{C'})^{\Z/h}$ and is stable under the multiplication with $\OO(C-p)$. 
Then our Manin triple still comes from a symmetric spherical sheaf of Lie algebras $\LL$ on $C$, 
equipped with a trivialization $\hat{L}_p\simeq \fg\ot\hat{\OO}$.
\end{prop}

\Pf . (i) Set $A=\OO(C-p)$.

\medskip

\noindent
{\bf Step 1. Recovering the curve as $\Proj$, and equipping it with a coherent sheaf $\LL$}.
Let us consider on $\fg\ot K$ the increasing filtration $F_n$ by the order of pole. Then we have the induced filtration $F_n\cap L$,
and the $A$-module structure on $L$ is compatible with this filtration and with the pole order filtration $(F_nA)$ on $A$.
Let us consider the graded module $\bigoplus_n F_n\cap L$ over $\RR(A):=\bigoplus_n F_nA$. It is well known 
that $\Proj \RR(A)$ is naturally isomorphic to $C$. 
Thus, from our graded module we get a coherent sheaf $\LL$ on $C$, which is a subsheaf of the constant sheaf $\fg\ot K_C$. 

\medskip

\noindent
{\bf Step 2. Compatibility of filtrations}. 
We have a natural isomorphism $L\to \LL(C-p)$ (compatible with embeddings into $\fg\ot K_C$) sending $F_n\cap L$ to $H^0(C,\LL(np))$.
Note that $F_0\cap L=L\cap \fg\ot\hat{\OO}=0$, while for every $n>0$ the natural map
$$(F_n\cap L)/(F_{n-1}\cap L)\to \fg\ot \OO(np)|_p$$
is an isomorphism.
We have a commutative diagram
\begin{diagram}
(F_n\cap L)/(F_{n-1}\cap L)&\rTo{\sim}& \fg\ot \OO(np)|_p\\
\dTo{}&&\dTo{}\\
H^0(C,\LL(np))/H^0(C,\LL((n-1)p))&\rTo{}&\LL(np)|_p&
\end{diagram}
which shows that the left vertical arrow
is injective. Hence, for every $n\ge 0$, $F_n\cap L=H^0(C,\LL(np))$.
Thus, the isomorphism $L\simeq \LL(C-p)$ is compatible with filtrations.

\medskip

\noindent
{\bf Step 3. Computing the cohomology of $\LL$}. 
We claim that the embedding $L\sub \fg\ot K$ induces an isomorphism $L\ot_A K\simeq \fg\ot K$.
Indeed, set $L_K=K\cdot L\sub \fg\ot K$ (so $L_K\simeq L\ot_A K$), and let $t$ be a formal parameter at $p$. Since $\fg\ot K=L_K+\fg\ot \hat{\OO}$, we get
$$\fg\ot t^{-1}\hat{\OO}\sub L_K+\fg\ot \hat{\OO},$$
or equivalently, $\fg\ot \hat{\OO}\sub L_K\cap (\fg\ot\hat{\OO})+\fg\ot t\hat{\OO}$. By Nakayama Lemma, it follows
that $\fg\ot \hat{\OO}=L_K\cap (\fg\ot\hat{\OO})$, Hence, $\fg\ot\hat{\OO}\sub L_K$, which implies our claim.

Therefore, $\dim_{K_C} L\ot_A K_C=\dim \fg$, so the embedding
$\LL\to \fg\ot K_C$ induces an isomorphism on stalks at the generic point,
$$\LL(C-p)\ot_A K_C\simeq \fg\ot K_C,$$
compatible with filtrations. 
On the other hand, we always have an isomorphism
$$\LL(C-p)\ot_A K_C\simeq \LL_p\ot_{\OO_{C,p}} K_C,$$
compatible with filtrations. Hence, we get an isomorphism
$$\LL_p\ot_{\OO_{C,p}} K_C\simeq \fg\ot K_C,$$ 
which is compatible with the pole/zero order filtration. Thus, it induces isomorphisms
$$\LL_p\ot_{\OO_{C,p}} K\simeq \fg\ot K, \ \ \LL_p\ot_{\OO_{C,p}}\hat{\OO}\simeq \fg\ot \hat{\OO}.$$
This implies that the complex $L\to (\fg\ot K)/(\fg\ot \hat{\OO})$ computes the cohomology of $\LL$,
so $H^*(C,\LL)=0$.

\medskip

\noindent
{\bf Step 4. Constructing the bracket and the $\om_C$-valued pairing on $\LL$}.
Since the Lie bracket is compatible with filtrations, it induces an $\OO$-linear Lie algebra structure on $\LL$. 
Next, we claim that the restriction to $L$ of the $K$-bilinear pairing
$$(\cdot,\cdot)_K: (\fg\ot K) \times (\fg\ot K)\to \om_{C,p}\ot K: (X\ot f, Y\ot g)\mapsto (X,Y)\cdot \eta\ot fg$$
takes values in $\om_C(C-p)$.
Indeed, since $L\sub \fg\ot K_C$, the induced pairing on $L$ takes values in rational $1$-forms on $C$.
Now we use the fact that $L$ is closed under multiplication with $\OO(C-p)$ and is isotropic with respect to $\Res_p(\cdot,\cdot)_K$,
so for any $l_1,l_2\in L$ and any $f\in \OO(C-p)$ we have
$$\Res_p(f\cdot (l_1,l_2)_K)=\Res_p(fl_1,l_2)_K=0.$$
By Lemma \ref{sing-om-lem}, this implies that $(l_1,l_2)_K\in \om_C(C-p)$.

Furthermore, the induced $\OO(C-p)$-bilinear pairing
$$L\ot L\to \om_C(C-p)$$
is compatible with the pole filtrations, so it induces a regular pairing
$$\kappa:\LL\ot \LL\to \om_C.$$
Since it is nondegenerate at $p$, by Lemma \ref{sph-Lie-def-lem}, $\LL$ is symmetric spherical.

\noindent
(ii) The proof is almost exactly the same as in (i). The main observation is that in all the arguments of (i)
we can use $\wt{\fg}(K_C):=(\fg\ot K_{C'})^{\Z/h}$ instead of $\fg\ot K_C$
(in fact, $\wt{\fg}(K_C)$ is a twisted form of $\fg\ot K_C$ trivialized on the extension $K_C\sub K_{C'}$),
where we equip $\wt{\fg}(K_C)$ with the pole/zero order filtration using the embedding \eqref{twisted-Lie-alg-K-emb}. 
Note that the embedding $i_q:K_{C'}\to K$ induces a pole/zero order filtration on $K_{C'}$ such that the completion gives $K$.
Hence, the embedding \eqref{twisted-Lie-alg-K-emb} also induces an isomorphism of the completion of $\wt{\fg}(K_C)$ with $\fg\ot K$.
In other words,
$$\wt{\fg}(K_C)\ot_{K_C} K\simeq \fg\ot K.$$
Now analogs of Steps 1-3 give a coherent sheaf $\LL$ (a subsheaf of the constant sheaf $\wt{\fg}(K_C)$), together with
isomorphisms 
$$L\simeq \LL(C-p), \ \  \LL_p\ot_{\OO_{C,p}} K_C\simeq \wt{\fg}(K_C),$$
compatible with filtrations, which implies that $H^*(\LL)=0$.

For the analog of Step 4, the only additional fact we need to know is that the restriction of the pairing $(\cdot,\cdot)_K$ to 
$\wt{\fg}(K_C)$ takes values in rational $1$-forms on $C$. To prove this let us write an explicit formula for the embedding
\eqref{twisted-Lie-alg-K-emb}. Note that all $\Z/h$-invariant elements in $\fg\ot K_{C'}$ are linear combinations of
the elements of the form $\sum_i \a^iX\ot \zeta^i(f)$, for $X\in \fg$, $f\in K_{C'}$.
The identification of $\fg\ot K$ with the subspace of $\Z/h$-invariants in $\fg\ot K[\Z/h]$ is given by
$$X\ot \varphi\mapsto \sum_j \a^jX\ot \varphi\cdot\zeta^j.$$
Using this it is easy to see that the image of $\wt{\fg}(K_C)$ in $\fg\ot K$ under \eqref{twisted-Lie-alg-K-emb} is spanned
by the elements of the form
$$\sum_i \a^iX\ot \iota_q(\zeta^i(f)).$$
Now we can calculate the pairing of two such elements
\begin{align*}
&(\sum_i \a^iX_1\ot \iota_q(\zeta^i(f_1)),\sum_i \a^iX_2\ot \iota_q(\zeta^i(f_2)))_K=\sum_{i,j}(\a^iX_1,\a^jX_2)\cdot
\eta\ot \iota_q(\zeta^i(f_1)\cdot\zeta^j(f_2))=\\
&\sum_m (\a^mX_1,X_2)\cdot \eta\ot (\sum_j \iota_q(\zeta^j(\zeta^m(f_1)\cdot f_2)).
\end{align*}
It remains to note that $\sum_j \zeta^j(\zeta^m(f_1)\cdot f_2)$ lies in $K_C\sub K_{C'}$.
\ed


\subsection{Belavin-Drinfeld's classification}\label{BD-class-sec}

Here we work over $\C$.

It is shown in \cite{BD-genYB} that every nondegenerate solution of the CYBE \eqref{CYBE-eq} (meromorphic on $U\times U$ for some
domain $U\sub \C$) is equivalent to a solution $r(x_1,x_2)$ that depends only on the difference $x_1-x_2$, so below we consider only such
solutions.

Let us recall the formula for nondegenerate trigonometric solutions of the CYBE from \cite{BD}.
Let $\fg$ be a simple Lie algebra equipped with a Coxeter automorphism $A$ of order $h$, 
We have the corresponding $\Z/h\Z$-grading $\fg=\bigoplus_{j\in \Z/h\Z} \fg_j$, where
$$\fg_j=\{x\in \fg \ |\ Ax=e^{j\cdot\frac{2\pi i}{h}}x\}.$$
We also set $\fh=\fg_0$.

Let $\Ga\sub \fh^*$ be the corresponding set of simple weights. A {\it Belavin-Drinfeld triple} (BD-triple) 
$(\Ga_1,\Ga_2,\tau)$ consists of two subsets $\Ga_1,\Ga_2\sub \Ga$ and a bijection
$\tau:\Ga_1\to \Ga_2$ preserving the scalar products. In addition, it is required that for any $\a\in \Ga$ the
expression $\tau^m(\a)$ is not defined for sufficiently large $m$.

Let $\Om_j\in \fg_j\ot \fg_{-j}$ denote the component of $\Om\in \fg\ot \fg$ in $\fg_j\ot \fg_{-j}$.
Note that $\Om^{21}_j=\Om_{-j}$.

The trigonometric $r$-matrix depends also on a continuous parameter $r_0\in \fh\ot \fh$ such that
$$r_0+r_0^{21}=\Om_0$$
$$(\tau\a\ot 1)(r_0)+(1\ot\a)(r)=0, \ \ \a\in\Ga_1.$$

Let us extend the operator $\tau$ to a bijection $\tau:\Pi_1\to \Pi_2$,
where for $i=1,2$, we denote by $\Pi_i\sub\Ga$ the set of all weights of the subalgebra generated by $e_\a$ with $\a\in\Ga_i$.
Note that $\Pi_i\cap (-\Pi_i)=\emptyset$ for $i=1,2$.
Using this extended $\tau$, the operator $\psi:\fg\to \fg$ is defined by
$$\psi(e_{\a})=\sum_{m\ge 1}e_{\tau^m\a}.$$

The result of Belavin-Drinfeld classification is that
\begin{equation}\label{BD-r-formula}
r(z)=r_0+\frac{1}{e^z-1}\sum_{j=0}^{h-1}\e^{jz/h}\Om_j-\sum_{j=1}^{h-1}e^{jz/h}(\psi\ot \id)\Om_j+
\sum_{j=1}^{h-1}e^{-jz/h}(\id\ot\psi)\Om_{-j}
\end{equation}
is a unitary solution of the CYBE, and in this way one gets all trigonometric solutions up to equivalence.

\subsection{Spherical sheaves of Lie algebras over the irreducible nodal curve of arithmetic genus $1$}\label{sph-Lie-nodal-sec}

Now, let $C$ be the irreducible nodal curve of arithmetic genus $1$ with the normalization $\pi:\P^1\to C$,
where $\pi(0)=\pi(\infty)$. We are going to construct a symmetric spherical sheaf of Lie algebra on $C$ corresponding
to each nondegenerate trigonometric solutions of the CYBE. 

Let $(\fg,A)$ be a simple Lie algebra with an automorphism of finite order, and let $h$ be the order of $A$.
Let $C'$ denote another copy of $C$, 
and let $\wt{C}\to C$ and $\wt{C'}\to C'$ be the normalizations.
We denote by $x$ the affine coordinate on $\wt{C}\setminus\{\infty\}$ and
by $y$ the coordinate on $\wt{C'}\setminus\{\infty\}$.

Let us consider the morphism
$$f:\wt{C'}\to C: y\mapsto y^h$$
(that factors through $\wt{C}$), 
and let $\zeta:\OO_{\wt{C'}}\to \OO_{\wt{C'}}$ be the automorphism $y\mapsto e^{-2\pi i/h}y$.
Let us define the action of $\Z/h$ on $\fg\ot f_*\OO_{\wt{C'}}$, so that the generator acts
by $A\ot \zeta$. We have a decomposition 
$$f_*\OO_{\wt{C'}}=\bigoplus_{j\in \Z/h} (f_*\OO_{\wt{C'}})_j,$$
where $\zeta$ acts as $e^{2\pi i j/h}$ on the $j$th summand. In particular, $(f_*\OO_{\wt{C'}})_j=\OO_{\wt{C}}$, where we
view $\OO_{\wt{C}}$ as a subsheaf of $f_*\OO_{\wt{C'}}$ due to the fact that $f$ factors through $\wt{C}$.
We also consider the points $0$ and $\infty$ on $\wt{C'}$ and the corresponding divisors $(0)$ and $(\infty)$,
and for any $m\in \Z$ we consider sheaves of $\OO_{\wt{C}}$-modules
$(f_*\OO_{\wt{C'}}(m(0)+m(\infty)))_{-j}$. 




We use the form $\eta_0=\frac{dx}{x}$ to define a trivialization of $\om_C$. 
As a smooth point $p\in C$ we take $x=1$.
We use a formal parameter $z$ at $p$, such that $e^z=x$. Note that $\eta_0=dz$.
We will always identify $\OO(C-p)$ with a subring of $\C(\!(z)\!)$ using this parameter.
Note that $\OO(C-p)$ is generated over $\C$ by $\frac{x}{(x-1)^2}$ and $\frac{x}{(x-1)^3}$.


Below we will consider divisors $m_0(0)+m_\infty(\infty)$ supported on $\{0,\infty\}\sub \wt{C'}$ and denote by
$\OO_{\wt{C'}}(m_0(0)+m_\infty(\infty))$ the corresponding subsheaves in the push-forward of $\OO_{\wt{C'}\setminus \{0,\infty\}}$.

\begin{thm}\label{main-Lie-thm}
Let $r(z)$ be the trigonometric $r$-matrix \eqref{BD-r-formula} associated with a Belavin-Drinfeld triple $(\Ga_1,\Ga_2,\tau)$ and with a continuous parameter $r_0$.
Then the corresponding Lie subalgebra $\fg(r)\sub \fg(\!(z)\!)$ is closed under the multiplication with $\OO(C-p)\sub \C(\!(z)\!)$.
Hence, the corresponding Manin triple comes from a symmetric spherical sheaf of Lie algebras $\LL$ on $C$ (see Proposition \ref{triple-to-sheaf-prop}). 
Furthermore, we have a $\Z/h$-grading $\LL=\bigoplus_{j\in\Z/h} \LL_j$, and 
we can realize $\LL$ as a subsheaf of $\Z/h$-graded Lie algebras in $f_*\OO_{\wt{C'}\setminus\{0,\infty\}}\ot \fg$, 
so that
$$\II\ot \fg_0\sub \LL_0\sub \OO_{\wt{C}}\ot \fg_0,$$
$$(f_*\OO_{\wt{C'}}(-(h+1)(0)-(h+1)(\infty)))_{-j}\ot \fg_j\sub \LL_j\sub (f_*\OO_{\wt{C'}}((h-1)(0)+(h-1)(\infty)))_{-j}\ot \fg_j $$
for $j\not\equiv 0$, where $\II\sub \OO_C$ is the ideal of the node.
\end{thm}

Recall that for the construction of $\fg(r)$ we introduce the second copy of the formal parameter $z$ which we call $t$.
We view $K_C$, as the subfield $\C(x)\sub \C(\!(t)\!)$, where $x=e^t$. 
Similarly, we view $K_{C'}$ as the subfield $\C(y)\sub \C(\!(t)\!)$, where $y=e^{t/h}$.
It is clear from \eqref{BD-r-formula} that $\fg(r)\sub \fg\ot K_{C'}$.

We will use the expansion
$$\frac{1}{e^{t-u}-1}=\frac{1}{xe^{-u}-1}=\sum_{m\ge 0}a_mu^m,$$
where $a_m\in \C(x)\sub \C(\!(t)\!)$. For example,
$$a_0=\frac{1}{x-1}, \ \ a_1=\frac{x}{(x-1)^2}, \ \ a_2=\frac{x}{(x-1)^3}.$$

Note that although the points $0$ and $\infty$ on $\wt{C'}$ appear symmetrically in the statement of Theorem \ref{main-Lie-thm},
we break this symmetry by using the coordinate $y$ centered at $0$, so it is not surprizing that in some elements of our proof the points $0$ and $\infty$ play slightly different roles. For example, this is so in the following technical result.

\begin{lem}\label{node-ideals-lem}
Let $I\sub \OO(C-p)$ denote the ideal of the node, and let 
$$I_2=H^0(\wt{C}-p,\OO(-(0)-2(\infty))).$$
Then we have
$$\OO(\wt{C}-p)=\operatorname{span}(1, (a_m)_{m\ge 0}),$$
$$I=\operatorname{span}(a_m \ |\ m\ge 1),$$
$$I_2=\operatorname{span}(a_m-\frac{1}{m!}a_1 \ |\ m\ge 2).$$
\end{lem} 

\Pf . Set $b_m=x^m/(x-1)^{m+1}$.
It is clear that $(b_m)_{m\ge 1}$ is a basis of $I$, whereas $(1,b_0,b_1,\ldots)$ is a basis of $\OO(\wt{C}-p)$.
Now the expansion 
$$\frac{1}{xe^{-u}-1}=\frac{1}{x-1}(1-\frac{x}{x-1}(1-e^{-u}))^{-1}=\sum_{m\ge 0} b_m(1-e^{-u})^m$$
shows that the $(a_m)_{m\ge 1}$ are expressed in terms of $(b_m)_{m\ge 1}$ via some invertible upper-triangular matrix.
This implies the statements involving $\OO(\wt{C}-p)$ and $I$.


Now let us check the statement about $I_2$.
Since $\operatorname{span}(a_m-\frac{1}{m!}a_1 \ |\ m\ge 2)$ has codimension $1$ in $I$, it is enough to check the inclusion
$a_m-\frac{1}{m!}a_1\in I_2$ for each $m\ge 2$.
To this end we note that $\sum_{m\ge 1}(a_m-\frac{1}{m!}a_1)u^m$ is the nonconstant part of the expansion in $u$ of
$$\frac{1}{xe^{-u}-1}-a_1(e^u-1)=\frac{1}{x-1}\cdot F(x,u),$$
where 
$$F(x,u)=(1-\frac{x}{x-1}(1-e^{-u}))^{-1}-\frac{x}{x-1}(e^u-1).$$
Thus, it is enough to check that all nonconstant terms of the expansion of $F(x,u)$ in $u$ vanish at $\infty$. In other words, we
need to check that $F(\infty,u)$ is constant. But
$$F(\infty,u)=e^u-(e^u-1)=1.$$
\ed

\medskip

\noindent
{\it Proof of Theorem \ref{main-Lie-thm}.}
Let us calculate $\fg(r)$. Let us denote by $r^{21}_{-j}$ the component of $r^{21}$ in $\fg_{-j}\ot \fg_j$. 
We have
$$r^{21}(t-u)_0=r_0^{21}+\frac{1}{xe^{-u}-1}\cdot\Om_0,$$
and for $1\le j\le h-1$,
$$r^{21}(t-u)_{-j}\cdot e^{\frac{ju}{h}}=\frac{y^j}{xe^{-u}-1}\cdot \Om_{-j}-y^j(\id\ot\psi)\Om_{-j}+y^{j-h}e^u(\psi\ot\id)\Om_{-j}.$$
Now we need to take the terms of expansion in $u$ and contract with a basis of $\fg^*$ in the first tensor component.
From this we immediately see that
$$\fg(r)=\bigoplus_{j\in \Z/h} \fg(r)_j, \text{ with } \fg(r)_j\sub y^jK_C\ot \fg_j.$$
Furthermore, we see that $\fg(r)_0$ is spanned by 
$$(e_i^*\ot \id)(r_0^{21})+a_0\ot e_i, \ \ (a_m\ot e_i)_{m\ge 1},$$
where $(e_i)$ is an orthonormal basis of $\fg_0$. On the other hand, 
for $1\le j\le h-1$, $\fg(r)_j$ is spanned by the elements
$$c^j_0(\alpha):=y^ja_0\ot e_{\a}-y^j\ot (e_{\tau \a}+e_{\tau^2\a}+\ldots)+y^jx^{-1}\ot (e_{-\tau^{-1}(-\a)}+e_{-\tau^{-2}(-\a)}+\ldots),$$
$$c^j_m(\alpha):=y^ja_m\ot e_{\a}+\frac{1}{m!}y^jx^{-1}\ot (e_{-\tau^{-1}(-\a)}+e_{-\tau^{-2}(-\a)}+\ldots), \text{ for } m\ge 1,$$
where $\a$ runs over all roots of $\fg_j$.

By Proposition \ref{triple-to-sheaf-prop}(ii) we just need to check that each $\fg(r)_j$ is closed under mutliplication by $\OO(C-p)$.
For this, it is enough to check closure under multiplication by $x/(x-1)^2$ and $x/(x-1)^3$.

Since $\operatorname{span}(a_m \ |\ m\ge 1)=I$, while $\OO(\wt{C}-p)=\operatorname{span}(1, (a_m)_{m\ge 0})$, 
we have 
$$\OO(C-p)\ot \fg\sub \fg(r)_0\sub \OO(\wt{C}-p)\ot \fg.$$
Since $I\cdot \OO(\wt{C}-p)\sub \OO(C-p)$, this implies that $\fg(r)_0$ is stable under multiplication by $I$. 

Now let us fix $j$, $1\le j\le h-1$.
Note that for $m\ge 2$,
$$c^j_m(\a)-\frac{1}{m!}c^j_1(\a)=y^j(a_m-\frac{1}{m!}a_1)\ot e_\a.$$
Hence, we get the inclusion
\begin{equation}\label{yj-I2-inclusion}
y^j I_2\ot \fg\sub  \fg(r)_j.
\end{equation}


Next, we claim that the following inclusions hold whenever the left-hand side is well defined:
\begin{equation}\label{yj-I-tau-a-inclusion}
y^j I\ot e_{\tau(\a)}\sub \fg(r)_j,
\end{equation}
\begin{equation}\label{yj-a1-a0-inclusion}
y^j(a_1-a_0)\ot e_{-\tau^{-1}(-\a)}\sub \fg(r)_j.
\end{equation}
Indeed, the first inclusion follows from the fact that for $m\ge 1$, one has
$$c^j_m(\tau(\a))=y^ja_m\ot e_{\tau(\a)}$$
since $\tau^{-1}$ is not defined on $-\tau(\a)$.
For the second inclusion we use the formula
\begin{equation}\label{c1-c0-difference}
c^j_1(\b)-c^j_0(\b)=y^j((a_1-a_0)\ot e_\b+1\ot e_{\tau(\b)}+\ldots).
\end{equation}
Applying this for $\b=-\tau^{-1}(-\a)$ and using the fact that $\tau$ is not defined on $-\tau^{-1}(-\a)$, we get
$$c^j_1(-\tau^{-1}(-\a))-c^j_0(-\tau^{-1}(-\a))=y^j(a_1-a_0)\ot e_{-\tau^{-1}(-\a)}$$
which proves our claim.

Now we are ready to check that $\fg(r)_j$ is stable under multiplication with $I$. 
We have to check the inclusion $I\cdot c^j_m(\a)\sub \fg(r)_j$ for every $\a$ and $m\ge 0$.
Since for $m\ge 2$, we have 
$$c^j_m(\a)-\frac{1}{m!}c^j_1(\a)\in I\ot \fg_j,$$
and since $I\cdot I\sub I_2$, the inclusion \eqref{yj-I2-inclusion} shows that it is enough to check the inclusions
$$I\cdot c^j_0(\a)\sub \fg(r)_j, I\cdot (c^j_1(\a)-c^j_0(\a))\sub \fg(r)_j.$$
Since 
$$I\cdot (a_1-a_0)=I\cdot \frac{1}{(x-1)^2}\sub I_2,$$
formula \eqref{c1-c0-difference} shows that for $f\in I$,
$$f\cdot c^j_1(\a)-c^j_0(\a)\equiv y^jf\ot (e_{\tau(\a)}+e_{\tau^2(\a)}+\ldots) \mod I_2\ot \fg_j.$$
Hence, the fact that this lies in $\fg(r)_j$ follows from \eqref{yj-I-tau-a-inclusion}.

Finally, since $a_0\cdot I\sub I_2$, using \eqref{yj-I-tau-a-inclusion}, we see that for $f\in I$
$$f\cdot c^j_0(\a)\equiv y^jx^{-1}f\ot (e_{-\tau^{-1}(-\a)}+e_{-\tau^{-2}(-\a)}+\ldots) \mod I_2\ot \fg_j.$$
It is enough to check this for $f=\frac{x}{(x-1)^2}$ and $f=\frac{x^2}{(x-1)^3}$.
In the latter case we have $x^{-1}f\in I_2$, so we are done. In the former case we get
$$\frac{x}{(x-1)^2}\cdot c^j_0(\a)\equiv y^j\frac{1}{(x-1)^2}\ot (e_{-\tau^{-1}(-\a)}+e_{-\tau^{-2}(-\a)}+\ldots) \mod I_2\ot \fg_j.$$
Since $\frac{1}{(x-1)^2}=a_1-a_0$, by \eqref{yj-a1-a0-inclusion}, this lies in $\fg(r)_j$.
\ed

\begin{ex}
In the case $\Ga_1=\Ga_2=\emptyset$ and $r_0=t_0/2$, we get from the above calculation,
$$\LL_j=(f_*\OO_{\wt{C'}})_{-j}\ot \fg_j \ \text{ for } j\not\equiv 0,$$
$$\LL_0=\FF\ot \fg_0,$$
where $\FF\sub \OO_{\wt{C}}$ is the $\OO_C$-submodule 
$$\FF:=\{f\in \OO_{\wt{C}} \ |\ f(0)+f(\infty)=0\}.$$
Note that the natural pairing $\LL\ot \LL\to \OO_{\wt{C}}$ factors through $\OO_C$,
since 
$$(f_*\OO_{\wt{C'}})_{-j}\cdot (f_*\OO_{\wt{C'}})_j\sub \II\sub \OO_C, \ \ \FF\cdot \FF\sub \OO_C.$$
This gives the pairing with values in $\om_C\simeq \OO_C$.

In the case when $A$ is an inner automorphism the solution corresponding to $\Ga_1=\Ga_2=\emptyset$ and $r_0=t_0/2$ is
equivalent to the standard solution 
$$r(z)=\frac{r^{21}e^z+r}{e^z-1}$$
coming from the standard quasitriangular structure $r$ on $\fg$. The same equivalence gives an isomorphism of our sheaf of algebras $\LL$ 
with the sheaf 
$$\LL_{st}\sub \fg\ot \OO_{\wt{C}}, \ \LL_{st}:=\{X\in \fg\ot \OO_{\wt{C}} \ |\ X(0)\in \fb_+, \ X(\infty)\in \fb_0, \ X(0)+X(\infty)\in \fn_-\oplus \fn_+\}$$
(see \ \cite[Prop.\ 2.13]{EK}).
\end{ex}




\end{document}